\documentclass[12pt]{article}
\usepackage{graphicx}
\usepackage{cite}
\usepackage[table]{xcolor}
\usepackage{booktabs}
\usepackage{multirow}
\usepackage{amsmath,amssymb}
\usepackage{pstricks}
\usepackage{amsmath,amssymb,subfigure,epsfig,bbm,stmaryrd,MnSymbol,mathrsfs}
\usepackage{subfigure}
\DeclareMathOperator{\sign}{sign}
\DeclareMathOperator{\supp}{supp}
\DeclareMathAlphabet{\mathpzc}{OT1}{pzc}{m}{it}

\begin{document}
\newcommand{\BigFig}[1]{\parbox{12pt}{\Huge #1}}
\newcommand{\BigZero}{\BigFig{0}}
\title{\bf Parametric Level Set Methods for Inverse Problems}   
\author{Alireza Aghasi, Misha Kilmer, Eric L. Miller\footnote{A. Aghasi and E. L. Miller are with the Department of Electrical Engineering, Tufts University. Misha Kilmer is with the Department of Mathematics, Tufts University, Medford, MA. Emails: \texttt{aaghas01@ece.tufts.edu, misha.kilmer@tufts.edu, elmiller@ece.tuts.edu}}\\
}         
\date{\today}    
\maketitle

\begin{abstract}
In this paper, a parametric level set method for reconstruction of
obstacles in general inverse problems is considered. General
evolution equations for the reconstruction of unknown obstacles
are derived in terms of the underlying level set parameters. We
show that using the appropriate form of parameterizing the level
set function results a significantly lower dimensional problem,
which bypasses many difficulties with traditional level set
methods, such as regularization, re-initialization and use of
signed distance function. Moreover, we show that from a
computational point of view, low order representation of the
problem paves the path for easier use of Newton and quasi-Newton
methods. Specifically for the purposes of this paper, we
parameterize the level set function in terms of adaptive compactly
supported radial basis functions, which used in the proposed
manner provides flexibility in presenting a larger class of shapes
with fewer terms. Also they provide a ``narrow-banding" advantage
which can further reduce the number of active unknowns at each
step of the evolution. The performance of the proposed approach is
examined in three examples of inverse problems, i.e., electrical
resistance tomography, X-ray computed tomography and diffuse
optical tomography.\\
\end{abstract}
{\scriptsize \bf DOI:\normalsize} 10.1137/100800208

\section{Introduction} \label{intro2}
Inverse problems arise in many applications of science and
engineering including e.g., geophysics
\cite{zhdanov2002geophysical,snieder1999inverse}, medical imaging
\cite{louis1992medical,webb2003introduction,arridge},
nondestructive evaluation
\cite{marklein2002linear,liu2003computational} and hydrology
\cite{yeh22review,carrera2005inverse,sun1994inverse}. In all cases
the fundamental problem is usually extracting the information
concerning the internal structure of a medium based on indirect
observations collected at the periphery, where the data and the
unknown are linked via a physical model of the sensing modality. A
fundamental challenge associated with many inverse problems is
ill-posedness (or ill-conditioning in the discrete case), meaning
that the solution to the problem is highly sensitive to noise in
the data or effects not captured by the physical model of the
sensor. This difficulty may arise due to the underlying physics of
the sensing system which in many cases (e.g., electrical impedance
tomography \cite{cheney1999electrical}, diffuse optical tomography
\cite{arridge}, inverse scattering \cite{dorn2006level}, etc)
causes the data to be inherently insensitive to fine scale
variations in the medium. This phenomenon makes such
characteristics difficult, if not impossible to recover stably
\cite{tikhonov1977solutions}. Another important factor causing the
ill-posedness is limitations in the distribution of the sensors
yielding sparse data sets that again do not support the recovery
of fine scale information uniformly in the region of interest
\cite{miller1996multiscale}. Many inverse problems of practical
interest in fact suffer from both of these problems. From a
practical point of view, left untreated, ill-posedness yields
reconstructions contaminated by high frequency, large amplitude
artifacts.

Coping with the ill-posedness is usually addressed through the use
of regularization
\cite{morozov1993regularization,engl1996regularization}. Based on
prior information about the unknowns, the regularization schemes
add constraints to the inverse problem to stabilize the
reconstructions. When the inverse problem is cast in a variational
framework, these regularization methods often take the form of
additive terms within the associated cost function and are
interpreted as penalties associated with undesirable
characteristics of the reconstruction. They may appear in various
forms such as imposing boundedness on the values of the unknown
quantity (e.g., Tikhonov or minimum norm regularizations
\cite{tikhonov1977solutions,golub2000tikhonov}) or penalizing the
complexity by adding smoothness terms (e.g., the total variation
regularization \cite{acar1994analysis}). These regularization
schemes are employed in cases where one seeks to use the data to
determine values for a collection of unknowns associated with a
dense discretization of the medium (e.g, pixels, voxels,
coefficients in a finite element representation of the unknown
\cite{rekanos1999inverse}).

For many problems, the fundamental objective of the process is the
identification and characterization of regions of interest in a
medium (tumors in the body \cite{bushberg2003essential},
contaminant pools in the earth \cite{james33optimal}, cracks in a
material sample \cite{stavroulakis2001inverse}, etc). For such
problems, an alternative to forming an image and then
post-processing to identify the region is to use the data to
directly estimate the geometry of the regions as well as the
contrast of the unknown in these regions. Problems tackled in this
way are known as the inverse obstacle or shape-based problems. For
earlier works in this area see \cite{kress1992numerical,
colton1996simple, kirsch1998characterization} and particularly the
more theoretical efforts by Kirsch \cite{kirsch1993domain} and
Kress {\it et al.} \cite{kress1994quasi}. Such processes usually
involve a rather simple parametrization of the shape and perform
the inversion based on using the domain derivatives mapping the
scattering obstacle to the observation. Relative to pixel-type
approaches, these geometric formulations tend to result in better
posed problems due to a more appropriate obstacle representation.
Moreover, in such problems the regularization can be either
performed implicitly through the parametrization or expressed in
terms of geometric constraints on the shape
\cite{kristensson1986inverse}. However, this class of shape
representation is not topologically flexible and the number of
components for the shape should be a priori known
\cite{santosa1996level}. They also create difficulties in
encountering holes and high curvature regions such as the corners.
These difficulties have lead over the past decade or so to the
development of shape-based inverse methods employing level
set-type representation of the unknowns.

The concept of level sets was first introduced by Osher and
Sethian in \cite{osher1988fronts}. This method was initially
designed for tracking the motion of a front whose speed depends on
the local curvature. The application of the level set approach to
inverse problems involving obstacles was discussed by Santosa in
\cite{santosa1996level}. One of the most attractive features of
the level set method is its ability to track the motion through
topological changes. More specifically, an a priori assumption
about the connectedness of the shapes of interest was no longer
required. Following Santosa's work, Litman {\it et al.} in
\cite{litman1998reconstruction} explored the reconstruction of the
cross-sectional contour of a cylindrical target. Two key
distinguishing points about this work were in the way that the
authors dealt with the deformation of the contour, and in their
use of the level set method to represent the contour. The shape
deformation method implemented in this work was enabled by a
velocity term, and lead to a closed-form derivative of the cost
functional with respect to a perturbation of the geometry. They
defined shape optimization as finding a geometry that minimizes
the error in the data fit. Later Dorn {\it et al.}~in
\cite{dorn2000shape} introduced a two-step shape reconstruction
method that was based on the adjoint field and level set methods.
The first step of this algorithm was used as an initialization
step for the second step and was mainly designed to deal with the
non-linearities present in the model. The second step of this
algorithm used a combination of the level set and the adjoint
field methods. Although inspired by the works of
\cite{osher1988fronts,santosa1996level,litman1998reconstruction},
the level set method used by Dorn {\it et al.}~was not based on a
Hamilton-Jacobi type equation, instead, an optimization approach
was employed, and an inversion routine was applied for solving the
optimization. The level set ideas in inverse problems were further
developed to tackle more advanced problems such as having shapes
with textures or multiple possible phases \cite{dorn2007level,
chan2004level}. Moreover, regarding the evolution of the level set
function where usually gradient descent methods are the main
minimization schemes applied, some authors such as Burger
\cite{burger2004levenberg} and Soleimani
\cite{soleimani2008computational} proposed using second order
convergent methods such as the Newton and quasi-Newton methods.

Although level set methods provide large degrees of flexibility in
shape representation, there are numerical concerns associated with
these methods. Gradient descent methods used in these problems
usually require a long evolution process. Although this problem
may be overcome using second order methods, the performance of
these methods for large problems such as 3D shape reconstructions
remains limited and usually gradient descent type methods are the
only option for such problems. Moreover re-initialization of the
level set function to keep it well behaved and velocity field
extensions to globally update the level set function through the
evolution are usually inevitable and add extra computational costs
and complexity to the problem \cite{osher2003level}. More detailed
reviews of level set methods in inverse problems can be found in
\cite{dorn2006level,burger2005survey}.

In all traditional level set methods already stated, the unknown
level set function belongs to an infinite dimensional function
space. From an implementation perspective, this requires the
discretization of the level set function onto a dense collection
of nodes. An alternative to this approach is to consider a finite
dimensional function space or a parametric form for the level set
function such as the space spanned by a set of basis functions.
Initially, Kilmer {\it et al.} in \cite{kilmer2004cortical}
proposed using a polynomial basis for this purpose in diffuse
optical tomography applications. In this approach, the level set
function is expressed in terms of a fixed order polynomial and
evolved through updating the polynomial coefficients at every
iteration. Parametrization of the level set function later
motivated some authors in field of mechanics to use it in
applications such as structural topology optimization
\cite{wang2006radial,wang2007extended,pingen2010parametric}. One
of the main contributions in this regard is the work by Wang {\it
et al.} in \cite{wang2006radial}. Here the level set function is
spanned by multiquadric radial basis functions as a typical basis
set in scattered data fitting and interpolation applications.
Authors showed that through this representation, the
Hamilton-Jacobi partial differential equation changes into a
system of ordinary differential equations and the updates for the
expansion coefficients may be obtained by numerically solving an
interpolation problem at every iteration.

More recently, the idea of parametric representation of the level
set function has been considered for image processing applications
\cite{gelas2007compactly,bernard2009variational}. Gelas {\it et
al.} in \cite{gelas2007compactly} used a similar approach as the
one by Wang {\it et al.} for image segmentation. As the basis set
they used compactly supported radial basis functions, which not
only reduce the dimensionality of the problem due to the
parametric representation, but also reduce the computation cost
through the sparsity that this class of functions provide. As
advantages of the method they showed that appropriate constraints
on the underlying parameters can avoid implementation of the usual
re-initialization. Also the smoothness of the solution is
guaranteed through the intrinsic smoothness of the underlying
basis functions, and in practice no further geometric constraints
need to be added to the cost function. As an alternative to this
approach, Bernard {\it et al.} in \cite{bernard2009variational}
parameterized the level set function in terms of B-splines. One of
the main advantages of their method was representing the the cost
function minimization directly in terms of the B-spline
coefficients and avoiding the evolution through the
Hamilton-Jacobi equation.

In this paper the general parametric level set approach for
inverse problems is considered. For an arbitrary parametrization
of the level set function, the evolution equation is derived in
terms of the underlying parameters. Since one of the main
advantages of this approach is low order representation of the
level set function, in practice the number of unknown parameters
in the problem is much less than the number of pixels (voxels) in
a traditional Hamilton-Jacobi type of level set method, therefore
we concentrate on faster converging optimization methods such as
the Newton or quasi-Newton type method. To represent the
parametric level set function, as in \cite{wang2006radial,
gelas2007compactly} we have proposed using radial basis functions.
However unlike the previous works employing radial basis functions
in a level set framework, in addition to the weighting
coefficients, the representation is adaptive with respect to the
centers and the scaling of the underlying radial basis functions.
This technique basically prevents using a large number of basis
terms in case that no prior information about the shape is
available. To fully benefit from this adaptivity, we further
narrow our choice by considering compactly supported radial basis
functions. Apparently, this choice would result sparsity in the
resulting matrices. However, we will discuss a behavior of these
functions which can be exploited to further reduce the number of
underlying basis terms and provide the potential to reconstruct
rather high curvature regions. The flexibility and performance of
proposed methods will be examined through illustrative examples.

The paper is structured as follows. In Section \ref{sec2} we
review shape-based inverse problems in a general variational
framework. Section \ref{sec3} is concerned with obtaining the
first and second order sensitivities of the cost function with
respect to the functions defining the shape. Based on the details
provided, a brief revision of the relevant traditional level set
methods is provided paving the path to use parametric level set
methods. In Section \ref{sec4} the general parametric level set
method will be introduced and evolution equations corresponding to
the underlying parameters are derived. In Section \ref{sec5} an
adaptive parametric level set based on the radial basis functions
will be proposed and the approach will be narrowed down to
compactly supported class of functions, due to their interesting
properties. Section \ref{sec6} will examine the method through
some examples in the context of electrical resistance tomography,
X-ray tomography and diffuse optical tomography. Finally in
Section \ref{sec7} some concluding remarks and suggestions will be
provided.

%

\section{Problem Formulation} \label{sec2}
\subsection{Forward Modelling}

The approach to modelling and determination of shape we consider
in this paper is rather general with the capability of being
applied to a range of inverse problems.  In Section \ref{sec6} we
specifically consider three applications: electrical resistance
tomography, limited view X ray tomography, and diffuse optical
tomography. The details of the specific problems will be provided
in the same section within the context of the examples themselves.
Up until that point, we have chosen to keep the discussion
general.

Consider $\Omega$ to be a compact domain in $\mathbb{R}^n$, $n\geq
2$, with Lipschitz boundary $\partial \Omega$. Further assume for
$\mathbf{x} \in \Omega$, a space dependent property $p(\mathbf{x})
\in \mathbb{S}_p$, where $\mathbb{S}_p$ is a Hilbert space. A
physical model $\mathcal{M}$ acts on a property of the medium
(e.g., electrical conductivity, mass density or optical
absorption), $p$, to generate an observation (or measurement)
vector $u$,
\begin{equation}
\label{eq1} u=\mathcal{M}(p),
\end{equation}
where $u$ itself belongs to some Hilbert space $\mathbb{S}_u$. In
most applications $u$ is a vector in $\mathbb{C}^k$, the space of
$k$ dimensional complex numbers where $k$ represents the number of
measurements and accordingly a canonical inner product is used.

As a convention throughout this paper, to keep the generality of
notation, the inner products and norms corresponding to any
Hilbert space $\mathbb{H}$, are subindexed with the notation of
the space itself, e.g. $\langle.,.\rangle_{\mathbb{H}}$.

\subsection{Inverse Problem}
The goal of an inverse problem is the recovery of information
about the property $p$ based on the data $u$. Here we consider a
variational approach where the estimate of $p$ is generated as the
solution to an optimization problem. The functional underlying the
problem is usually comprised of two terms. The first term demands
that the estimate of $p$ be consistent with the data in a
mathematically precise sense.  As is well known however, many
interesting inverse problems are quite ill-posed. This means when
the data consistency is our only concern, the resulting estimate
of $p$ could be quite far from the truth, corrupted by artifacts
such as high frequency and large amplitude oscillations. Hence, an
additional term (or terms) are required in the formulation of the
variation problem which capture our prior knowledge concerning the
expected behavior of the $p$ in $\Omega$. Such terms serve to
stabilize (or regularize) the inverse problem. Defining a residual
operator as
\begin{equation}
\label{eq2} \mathcal{R}(p)=\mathcal{M}(p)-u,
\end{equation}
the inverse problem is formulated in the following manner
\begin{equation}
\label{eq3} \underset{p}{\min}\: \ \mathscr{F}(p)=
\frac{1}{2}\|\mathcal{R}(p)\|_{\mathbb{S}_u}^2+\mathscr{L}\!(p),
\end{equation}
where $\mathscr{L}$ is the regularization functional. Appropriate
choice of $\mathscr{L}$ is usually based on properties of the
problem. In the typical case where the unknown property $p$ is
represented as a dense collection of pixels (or voxels) in
$\Omega$, the regularization penalties are used to enforce
smoothness and boundedness of $p$ in $\Omega$. These are usually
considered in the framework of Tikhonov and total variation
regularizations \cite{tikhonov1977solutions, acar1994analysis}. An
alternative approach that has been of great interest in recent
years is based on geometric parameterizations of the unknown. Here
the regularization penalties are either embedded in the nature of
the unknown or expressed as geometric constraints on the unknown.
No matter which approach is used, usually in defining
$\mathscr{L}$ different spaces and their corresponding norms may
be used. For a more detailed review of such methods an interested
reader is referred to \cite{vogel2002computational,dorn2006level}.
Clearly when the number of parameters involved in the problem is
sufficiently small, an underdetermined inverse problem can be made
overdetermined, and in this sense the problem becomes better
posed. Since the parametrization idea we will put forth in this
paper is empirically found to be well-posed enough that no
necessary regularization terms need to be added to the cost
function, $\mathscr{L}$ will be neglected in our future
discussions of $\mathscr{F}$.
\subsection{A Shape-Based Approach For the Unknown Parameter}
For a large class of ``shape-based" inverse problems
\cite{santosa1996level, kress1998inverse, dorn2006level,
miller2000new}, it is natural to view $p(\mathbf{x})$ as being
comprised of two classes, i.e., foreground and background. The
problem then amounts to determination of the boundary separating
these classes as well as characteristics of the property values in
each class. In the simplest case, $p(\mathbf{x})$ is piecewise
constant while in more sophisticated cases $p(\mathbf{x})$ may be
random and characterized by different probabilistic models in the
two regions.  In this paper, we assume that over each region
$p(\mathbf{x})$ is at least differentiable. The property of
interest in this case is usually formulated through the use of a
characteristic function. Given a closed domain $D\subseteq \Omega$
with corresponding boundary $\partial D$, the characteristic
function $\chi_D$ is defined as
\begin{equation}
\label{eq4} \chi_D(\mathbf{x}) = \left\{
\begin{array}{rl}
1 &  \mathbf{x} \in D\\
0 &  \mathbf{x} \in \Omega\setminus D .\\
\end{array} \right.
\end{equation}
Accordingly, the unknown property $p(\mathbf{x})$ can be defined
over the entire domain $\Omega$ as
\begin{equation}
\label{eq5}
p(\mathbf{x})=p_i(\mathbf{x})\chi_D(\mathbf{x})+p_o(\mathbf{x})\big(1-\chi_D(\mathbf{x})\big).
\end{equation}
Unlike $p(\mathbf{x})$ which is clearly not differentiable along
$\partial D$, $p_i(\mathbf{x})$ indicating the property values
inside $D$ and $p_o(\mathbf{x})$ denoting the values outside, are
assumed to be smooth functions, and as mentioned earlier, at least
belonging to $C^{1}(\Omega)$.

In a shape-based approach, finding $\partial D$ is a major
objective. As (\ref{eq5}) and (\ref{eq4}) show, $p(\mathbf{x})$ is
implicitly related to $D$ and to find $\partial D$, a more
explicit way of relating them should be considered. In this regard
the idea of using a level set function proves to be especially
useful \cite{santosa1996level}. Here $\partial D$ is represented
as some level set of a Lipschitz continuous function $\phi :
\Omega\rightarrow \mathbb{R}$. When the zero level set is
considered, $\phi(\mathbf{x})$ is related to $D$ and $\partial D$
via
\begin{equation}
\label{eq6} \left\{
\begin{array}{rl}
\phi(\mathbf{x})>0 &  \forall \mathbf{x} \in D \\
\phi(\mathbf{x})=0 &  \forall \mathbf{x} \in \partial D \\
\phi(\mathbf{x})<0 &  \forall \mathbf{x} \in \Omega\setminus D. \\
\end{array} \right.
\end{equation}
Making the use of a Heaviside function, defined as
$H(.)=\frac{1}{2}\big(1+\sign(.)\big)$, the function
$p(\mathbf{x})$ can be represented as
\begin{equation}
\label{eq7}
p(\mathbf{x})=p_i(\mathbf{x})H\big(\phi(\mathbf{x})\big)+p_o(\mathbf{x})\Big(1-H\big(\phi(\mathbf{x})\big)\Big).
\end{equation}
This equation in fact maps the space of unknown regions $D$ into
the space of unknown smooth functions $\phi$.
\section{Inversion as a Cost Function Minimization}\label{sec3}
In this section we develop the mathematical details of the
minimization problem (\ref{eq3}) when a shape-based approach as
(\ref{eq7}) is considered. Most current methods use the first and
second order sensitivities of the cost function with respect to
the unknown parameter to perform the minimization
\cite{burger2005survey}. Based on the general details provided, we
will briefly revisit the traditional level set approaches relevant
to this paper in Section \ref{revisit_sec}, since understanding
the details better justifies the use of parametric level set
methods.
\subsection{Cost Function Variations Due to the Unknowns}
We begin by assuming that the first and second order Fr\'{e}chet
derivatives of $\mathcal{R}(p)$ exist and denote them as
$\mathcal{R}'(p)[\hspace{.7mm}.\hspace{.7 mm}]$ and
$\mathcal{R}''(p)[.,.]$. The first order Fr\'{e}chet derivative of
a function (if it exists) is a bounded and linear operator. The
second order Fr\'{e}chet derivative is also bounded but bilinear,
which means the operator acts on two arguments and is linear with
respect to each \cite{berger1977nonlinearity}.

For an arbitrary variation $\delta p \in \mathbb{S}_p$ and the
real scalar $\varepsilon$, using the generalized Taylor expansion
we have
\begin{equation}
\label{eq8} \mathcal{R}(p+\varepsilon\delta
p)=\mathcal{R}(p)+\varepsilon\mathcal{R}'(p)[\delta p]
+\frac{\varepsilon^2}{2} \mathcal{R}''(p)[\delta p, \delta p] +
O(\varepsilon^3).
\end{equation}
Rewriting $\mathscr{F}(p)$ as
\begin{equation}
\label{eq9} \mathscr{F}(p)=\frac{1}{2}\langle
\mathcal{R}(p),\mathcal{R}(p)\rangle_{\mathbb{S}_u},
\end{equation}
and recalling the fact that $\langle u_1,u_2
\rangle_{\mathbb{S}_u}=\overline{\langle u_2,u_1
}\rangle_{\mathbb{S}_u}$ for $ u_1, u_2 \in \mathbb{S}_u$ and
overline denoting complex conjugate, the variations of the cost
function with respect to the variations of $p$ can be derived as
\begin{equation}
\label{eq10} \mathscr{F}(p+\varepsilon\delta
p)=\mathscr{F}(p)+\varepsilon\mathscr{F}'(p)[\delta p]
+\frac{\varepsilon^2}{2} \mathscr{F}''(p)[\delta p, \delta p] +
O(\varepsilon^3),
\end{equation}
where for $ p_1, p_2 \in \mathbb{S}_p$
\begin{equation}
\label{eq11} \mathscr{F}'(p)[p_1]=\mathscr{R}\!\mathpzc{e}\langle
\mathcal{R}'(p)[ p_1],\mathcal{R}(p)\rangle_{\mathbb{S}_u}
\end{equation}
and
\begin{equation}
\label{eq12}
\mathscr{F}''(p)[p_1,p_2]=\mathscr{R}\!\mathpzc{e}\langle
\mathcal{R}'(p)[p_1],\mathcal{R}'(p)[
p_2]\rangle_{\mathbb{S}_u}+\mathscr{R}\!\mathpzc{e}\langle
\mathcal{R}''(p)[ p_1, p_2],\mathcal{R}(p)\rangle_{\mathbb{S}_u}.
\end{equation}
The notation $\mathscr{R}\!\mathpzc{e}$ indicates the real part of
the corresponding quantity. Denoting $\mathcal{R}'(p)^*[\;.\;]$ as
the adjoint operator between $\mathbb{S}_u$ and $\mathbb{S}_p$ as
\begin{equation}\label{eq13}
\langle
\hat{u},\mathcal{R}'(p)[\hat{p}]\rangle_{\mathbb{S}_u}=\langle
\mathcal{R}'(p)^*[\hat{u}],\hat{p}\rangle_{\mathbb{S}_p}, \qquad
\forall \hat{u} \in \mathbb{S}_u ,\forall \hat{p} \in
\mathbb{S}_p,
\end{equation}
 (\ref{eq11}) can be written as
\begin{equation}
\label{eq14} \mathscr{F}'(p)[p_1]=\mathscr{R}\!\mathpzc{e}\langle
\mathcal{R}'(p)^* [\mathcal{R}(p)],p_1\rangle_{\mathbb{S}_p}.
\end{equation}
Equations (\ref{eq11}), (\ref{eq14}) and (\ref{eq12}) are in fact
the first and second order Fr\'{e}chet derivatives of
$\mathscr{F}$ with respect to $p$. In a more general context (and
indeed one which we shall use in Section \ref{sec4}), $p$ itself
can be the map over some variable $v$ from another Hilbert space
$\mathbb{S}_v$ into $\mathbb{S}_p$, i.e., $p(v): \mathbb{S}_v
\rightarrow \mathbb{S}_p$. Assuming the existence of the first and
second order Fr\'{e}chet derivatives of $p$ with respect to $v$,
denoted as $p'(v)[\hspace{.7mm}.\hspace{.7 mm}]$ and
$p''(v)[.,.]$, the first and second order Fr\'{e}chet derivatives
of $\mathscr{F}$ with respect to $v$ can be obtained using the
chain rule as
\begin{equation}
\label{eq16}
\mathscr{F}'(v)[v_1]=\mathscr{F}'(p)\big[p'(v)[v_1]\big]
\end{equation}
and
\begin{equation}
\label{eq17} \mathscr{F}''(v)[v_1,
v_2]=\mathscr{F}''(p)\big[p'(v)[v_1],p'(v)[v_2]\big]+\mathscr{F}'(p)\big[p''(v)[v_1,v_2]\big],
\end{equation}
where $ v_1, v_2 \in \mathbb{S}_v$. Equations (\ref{eq16}) and
(\ref{eq17}) themselves can be easily expressed in terms of
$\mathcal{R}(p)$ and its derivatives using (\ref{eq11}) and
(\ref{eq12}). These equations will be used later as the key
equations in finding the sensitivities in our parametric level set
representation of $p$.

\subsection{Pixel Based Minimizations (Revisiting Traditional Level Set
Methods)}\label{revisit_sec} In the specific context of the
shape-based inverse problems of interest here, the first and
second order sensitivities of the cost function $\mathscr{F}$ with
respect to the functions defining $p(\mathbf{x})$ in (\ref{eq7})
i.e., $\phi(\mathbf{x})$, $p_i(\mathbf{x})$ and $p_o(\mathbf{x})$,
can be used to form a minimization process. Based on the order of
the sensitivities available, first order optimization methods such
as gradient descent or second order methods such as Newton or
quasi-Newton techniques can be implemented.

For simplicity in reviewing the current methods, we assume that
$p_i$ and $p_o$ are known a priori and only the shape (i.e., the
zero level set of $\phi$) is unknown (see
\cite{villegas2006simultaneous, feng2003curve, dorn2007level} for
details on the recovery of both the shape as well as the contrast
function). In an evolution approach it is desired to initialize a
minimization process with some level set function $\phi_0$ and
evolve the function to find a $\phi$ which minimizes
$\mathscr{F}$. To take into account the concept of evolution, an
artificial time is defined where the level set function at every
time frame $t\geq 0$ is rewritten as $\phi(\mathbf{x};t)$ and the
zero level set of $\phi(\mathbf{x};t)$ is denoted as $\partial
D_t$. A straightforward differentiation of $\phi(\mathbf{x};t)=0$
with respect to $t$ yields to the Hamilton-Jacobi type equation
\begin{equation}
\label{eq18} \frac{\partial \phi}{\partial
t}+V(\mathbf{x};t)\cdot\nabla \phi =0
\end{equation}
for the points on $\partial D_t$ where
$V(\mathbf{x};t)=\mbox{d}\mathbf{x}/\mbox{d}t$. To move the
interface in the normal direction, $V(\mathbf{x};t)$ should be
chosen as $v(\mathbf{x};t)\vec{n}(\mathbf{x};t)$ where $v$ is a
scalar speed function and $\vec{n}=\nabla \phi/|\nabla \phi|$ is
the unit outward vector on $\partial D_t$. Incorporating this into
the minimization of $\mathscr{F}$, the speed function for the
points on $\partial D_t$, denoted as $\tilde{v}$, can be chosen to
be in the steepest descent direction of $\mathscr{F}$ which is
\cite{dorn2006level}
\begin{equation}
\label{eq19} \tilde{v}=-\mathscr{R}\!\mathpzc{e}\{
(p_o-p_i)\overline{\mathcal{R}'(p)^*
\left[\mathcal{R}(p)\right]}\}.
\end{equation}
As (\ref{eq19}) is only valid for $\mathbf{x} \in \partial D$, a
velocity extension should be performed to extend $\tilde{v}$ to
$v$ defined over the entire domain $\Omega$ and therefore capable
of globally evolve the level set function \cite{osher2003level}.
Beside this classical level set approach, other ways of
representing the speed functions and performing the minimization
process are proposed
\cite{miled2007projection,burger2004levenberg,
park2009reconstruction}. For example, Hadj Miled and Miller
\cite{miled2007projection} proposed a normalized version of the
classic speed function in the context of electrical resistance
tomography.

Some authors have also proposed using Newton type methods to
update $\tilde{v}$ at every iteration (e.g., see
\cite{santosa1996level,burger2004levenberg,soleimani2006level}).
Analogous to (\ref{eq2}) and (\ref{eq9}), in these problems the
residual operator and the cost function are usually written
directly as $\mathcal{R}(D)$ and $\mathscr{F}(D)$, functions of
the shape itself. Assuming the existence of the first and second
order shape derivatives \cite{novruzi2002structure}, denoted as
$\mathscr{F}'(D)[\hspace{.7mm}.\hspace{.7 mm}]$ and
$\mathscr{F}''(D)[.,.]$, at every time step the Newton update
$\tilde{v}$ is obtained by solving \cite{burger2005survey}
\begin{equation}
\label{eq20} \mathscr{F}''(D)[w,\tilde{v}]+\mathscr{F}'(D)[w]=0
\qquad \forall w \in \mathbb{S}_D.
\end{equation}
Here $\mathbb{S}_D$ is an appropriate Hilbert space such as
$L^2(\partial D)$, which may depend on the current shape
\cite{burger2004levenberg,burger2003framework}. Considering the
general forms of the derivatives as (\ref{eq11}) and (\ref{eq12}),
in a Gauss-Newton method the second derivatives of $\mathcal{R}$
are disregarded and (\ref{eq20}) becomes
\begin{equation}
\label{eq21} \mathscr{R}\!\mathpzc{e}\langle
\mathcal{R}'(D)[\tilde{v}],\mathcal{R}'(D)[
w]\rangle_{\mathbb{S}_u}+\mathscr{R}\!\mathpzc{e}\langle
\mathcal{R}'(D)[w],\mathcal{R}(D)\rangle_{\mathbb{S}_u}=0 \qquad
\forall w \in \mathbb{S}_D.
\end{equation}
Furthermore, to avoid ill conditioning, in a Levenberg-Marquardt
approach (\ref{eq21}) is regularized as
\begin{equation}
\label{eq22} \mathscr{R}\!\mathpzc{e}\langle
\mathcal{R}'(D)[\tilde{v}],\mathcal{R}'(D)[
w]\rangle_{\mathbb{S}_u}+\mathscr{R}\!\mathpzc{e}\langle
\mathcal{R}'(D)[w],\mathcal{R}(D)\rangle_{\mathbb{S}_u}+\lambda\langle
\tilde{v}, w \rangle_{\mathbb{S}_D}=0 \qquad \forall w \in
\mathbb{S}_D.
\end{equation}
in which for $\lambda>0$, the equation is shown to be well-posed
\cite{burger2004levenberg}. Similar to the previous approach, once
$\tilde{v}$ is obtained, a velocity extension is performed to
result a globally defined $v$ which can be used to update $\phi$.
However, although using second order methods can reduce the number
of iterations in finding a minima, compared to gradient descent
methods, they do not necessarily reduce the computation load.

From an implementation perspective there are some concerns using
the aforementioned methods. The gradient descent method usually
requires many iterations to converge and performances become poor
for low sensitivity problems \cite{miled2007projection}. Although
using Newton and quasi-Newton methods to update the level set
function increases the convergence rate, they are usually
computationally challenging and for large problems and relatively
finer grids, a large system of equations must be solved at every
iteration. Also for both types of methods, there are usually added
complications of the level set function re-initialization and
speed function extension. The approach that we will put forth in
the next section is capable of addressing these problems. It is
low order and numerically speaking, the number of unknowns
involved in the problem are usually much less than the number of
grid points and hence allows us to easily use second order
methods. Moreover, our proposed method does not require
re-initialization of the level set function, speed function
extension or even length-type regularization as a common
regularization in many shape-based methods (e.g., see
\cite{ito2001level,miled2007projection,dorn2006level}) and our
level set function remains well behaved through the corresponding
evolution process.

\section{A Parametric Level Set Approach}\label{sec4}

As discussed earlier, in most current shape-based methods
$\phi(\mathbf{x})$ is represented by function values on a dense
discretization of $\mathbf{x}$-space as part of a discretization
of the underlying evolution equation or Newton-type algorithm.
Consider now the level set function to be still a function of
$\mathbf{x}$ but also a function of a parameter vector
$\boldsymbol{\mu}=(\mu_1,\mu_2,\cdots,\mu_m) \in \mathbb{R}^m$. In
this case we define the continuous Lipschitz function $\phi :
\Omega \times \mathbb{R}^m\rightarrow \mathbb{R}$, as a
\emph{parametric level set} (PaLS) representation of $D$ if for a
$c \in \mathbb{R}$
\begin{equation}
\label{eq24} \left\{
\begin{array}{rl}
\phi(\mathbf{x},\boldsymbol{\mu})>c &  \forall \mathbf{x} \in D \\
\phi(\mathbf{x},\boldsymbol{\mu})=c &  \forall \mathbf{x} \in \partial D \\
\phi(\mathbf{x},\boldsymbol{\mu})<c &  \forall \mathbf{x} \in \Omega\setminus D. \\
\end{array} \right.
\end{equation}
In the PaLS approach we assume that the general form of
$\phi(\mathbf{x},\boldsymbol{\mu})$ is known and the specification
of $\boldsymbol{\mu}$ can explicitly define the level set function
over the entire domain $\Omega$. In other words the evolution of
$\phi$ required to solve the underlying inverse problems is
performed via the evolution of $\boldsymbol{\mu}$. An example of a
PaLS function is a basis expansion with known basis functions and
unknown weights and as will be shown later in this paper we
considerably expand on this notion. We call $\boldsymbol{\mu}$ the
PaLS parameters.

To setup the problem using the PaLS approach, consider momentarily
that $p_i(\mathbf{x})$ and $p_o(\mathbf{x})$ are known a priori
(later we will appropriately take away this restriction). Based on
(\ref{eq24}) $p$ is written as
\begin{equation}
\label{eq25}
p(\mathbf{x},\boldsymbol{\mu})=p_i(\mathbf{x})H\big(\phi(\mathbf{x},\boldsymbol{\mu})-c\big)+p_o(\mathbf{x})\Big(1-H\big(\phi(\mathbf{x},\boldsymbol{\mu})-c\big)\Big).
\end{equation}
Under this model, we now view $\mathscr{F}$ in (\ref{eq9}) as a
function of $\boldsymbol{\mu}$, i.e.
$\mathscr{F}(\boldsymbol{\mu}) : \mathbb{R}^m\rightarrow
\mathbb{R}$. Therefore unlike the classic level set approach, the
unknown is no longer the function $\phi$, but a vector belonging
to $\mathbb{R}^m$, where $m$ is usually much smaller than the
number of unknowns associated with a discretization of
$\mathbf{x}$-space. With this model and assuming (1) $\phi$ is
sufficiently smooth with respect to the elements of
$\boldsymbol{\mu}$ and (2) the discontinuous Heaviside function is
replaced by a $C^2$ approximation (e.g.,
\cite{zhao1996variational}) denoted as $H_{rg}$, we now proceed to
formulate a second order approach for the minimization of
$\mathscr{F}(\boldsymbol{\mu})$. To begin, rewriting (\ref{eq25})
with $H_{rg}$ and taking a derivative with respect to $\phi$
yields
\begin{equation}
\label{eq26} \frac{\partial p } {\partial
\phi}=(p_i-p_o)\delta_{rg}(\phi-c),
\end{equation}
where $\delta_{rg}(.)$ is accordingly the regularized version of
the Dirac delta function (see examples in
\cite{tian2008basis,zhao1996variational}). Using the chain rule
gives
\begin{equation}
\label{eq27} \frac{\partial p } {\partial \mu_j}=\frac{\partial p
} {\partial \phi}\frac{\partial \phi } {\partial
\mu_j}=(p_i-p_o)\delta_{rg}(\phi-c)\frac{\partial \phi } {\partial
\mu_j}.
\end{equation}
Now using (\ref{eq27}) with (\ref{eq16}) and (\ref{eq14}), the
gradient vector for $\mathscr{F}$ is
\begin{align}
\label{eq28} \nonumber \frac{\partial \mathscr{F}}{\partial
{\mu_j}}&=\mathscr{F}'(p)[\frac{\partial p}{\partial
\mu_j}]\\&=\mathscr{R}\!\mathpzc{e}\langle
\mathcal{R}(p),\mathcal{R}'(p)[\frac{\partial p}{\partial
\mu_j}]\rangle_{\mathbb{S}_u}\\&=\label{eq28.x}\mathscr{R}\!\mathpzc{e}\langle
\mathcal{R}'(p)^*
[\mathcal{R}(p)],(p_i-p_o)\delta_{rg}(\phi-c)\frac{\partial \phi }
{\partial \mu_j}\rangle_{\mathbb{S}_p}.
\end{align}
We denote by $\mathbf{J}_{\boldsymbol{\mu}}(\mathscr{F})$ the
gradient of $\mathscr{F}$ with respect to the parameter vector
$\boldsymbol{\mu}$. With this notation a gradient descent equation
can be formed to evolve the PaLS function as
\begin{equation}
\label{eq28a}
\boldsymbol{\mu}^{(t+1)}=\boldsymbol{\mu}^{(t)}-\lambda^{(t)}
\mathbf{J}_{\boldsymbol{\mu}}(\mathscr{F})\big
|_{\boldsymbol{\mu}=\boldsymbol{\mu}^{(t)}} \qquad t\geq 0,
\end{equation}
where $\lambda^{(t)}>0$ is the iteration step
\cite{bertsekas1999nonlinear}) and (\ref{eq28a}) is assumed to be
initialized with some $\boldsymbol{\mu^{(0)}}$. Although gradient
decent is relatively simple to implement, it is known to be slow
to converge and can suffer from difficulties associated with
scaling of the parameters \cite{dennis1996numerical}. Moreover,
the use of gradient decent fails to take advantage of one of the
primary benefits of the PaLS idea; namely the ability to specify a
level set function using a small (relative to a discretization of
$\mathbf{x}$-space) number of parameters. Under this model, it
becomes feasible and indeed useful to consider higher order
optimization methods such as Newton or quasi-Newton methods which
are faster in convergence and robust with respect to sensitivity
scalings of different parameters \cite{dennis1996numerical}. These
methods usually use the information in the Hessian, which we now
derive. To calculate the elements of the Hessian matrix for
$\mathscr{F}$ using (\ref{eq27}) we have
\begin{equation}
\label{eq29} \frac{\partial^2 p } {\partial \mu_j
\partial \mu_k}=(p_i-p_o)\big(\delta_{rg}(\phi-c)\frac{\partial^2 \phi }
{\partial \mu_j \partial \mu_k} +
\delta'_{rg}(\phi-c)\frac{\partial \phi } {\partial \mu_j}
\frac{\partial \phi}{\partial \mu_k} \big),
\end{equation}
where $\delta'_{rg}(.)$ is the derivative of the regularized Dirac
delta function. Based on (\ref{eq17}) and (\ref{eq12}) we have
\begin{align}
\label{eq30} \nonumber \frac{\partial^2 \mathscr{F}}{\partial
\mu_j\partial \mu_k}&= \mathscr{F}''(p)\big[\frac{\partial
p}{\partial \mu_j},\frac{\partial p}{\partial
\mu_k}\big]+\mathscr{F}'(p)\big[\frac{\partial^2 p } {\partial
\mu_j
\partial \mu_k}\big]
\\\nonumber &=\mathscr{R}\!\mathpzc{e}\langle
\mathcal{R}'(p)[ \frac{\partial p}{\partial \mu_j}],
\mathcal{R}'(p)[\frac{\partial p}{\partial
\mu_k}]\rangle_{\mathbb{S}_u}+\mathscr{R}\!\mathpzc{e}\langle
\mathcal{R}''(p)[ \frac{\partial p}{\partial \mu_j},
\frac{\partial p}{\partial
\mu_k}],\mathcal{R}(p)\rangle_{\mathbb{S}_u}\\&\hspace{5.2cm}+\mathscr{R}\!\mathpzc{e}\langle
\mathcal{R}'(p)^* [\mathcal{R}(p)],  \frac{\partial^2 p}{\partial
\mu_j\partial \mu_k} \rangle_{\mathbb{S}_p}.
\end{align}
This equation is in fact the exact expression for the elements of
the Hessian matrix for $\mathscr{F}$. However as mentioned
earlier, in methods such as the Gauss-Newton or
Levenberg-Marquardt, to reduce the computation cost the Hessian is
approximated in that the terms containing second order derivatives
are disregarded. Following that approach here, (\ref{eq30})
becomes
\begin{align}
\label{eq31}  \frac{\partial^2 \mathscr{F}}{\partial \mu_j\partial
\mu_k}&\simeq\mathscr{R}\!\mathpzc{e}\langle \mathcal{R}'(p)[
\frac{\partial p}{\partial \mu_j}], \mathcal{R}'(p)[\frac{\partial
p}{\partial
\mu_k}]\rangle_{\mathbb{S}_u}\\&=\label{eq31.b}\mathscr{R}\!\mathpzc{e}\langle
\mathcal{R}'(p)[ (p_i-p_o)\delta_{rg}(\phi-c)\frac{\partial \phi }
{\partial \mu_j} ],
\mathcal{R}'(p)[(p_i-p_o)\delta_{rg}(\phi-c)\frac{\partial \phi }
{\partial \mu_k}] \rangle_{\mathbb{S}_u}.
\end{align}
We denote as $\mathbf{\tilde{H}}_{\boldsymbol{\mu}}(\mathscr{F})$
the approximate Hessian matrix, the elements of which are obtained
through (\ref{eq31}). Having this in hand, a stable and faster
converging PaLS function evolution can be proposed by solving the
following Levenberg-Marquardt equation for
$\boldsymbol{\mu}^{(t+1)}$
\begin{equation}
\label{eq31c} \Big [
\mathbf{\tilde{H}}_{\boldsymbol{\mu}}(\mathscr{F})\big
|_{\boldsymbol{\mu}=\boldsymbol{\mu}^{(t)}}+\lambda^{(t)}
\mathbf{I}\Big](\boldsymbol{\mu}^{(t+1)}-\boldsymbol{\mu}^{(t)})=-\mathbf{J}_{\boldsymbol{\mu}}(\mathscr{F})\big
|_{\boldsymbol{\mu}=\boldsymbol{\mu}^{(t)}} \qquad t\geq 0.
\end{equation}
Here $\lambda^{(t)}$ is a small positive number chosen at every
iteration and $\mathbf{I}$ is the identity matrix. In fact
referring to (\ref{eq31}) we can see that the approximate Hessian
matrix $\mathbf{\tilde{H}}_{\boldsymbol{\mu}}(\mathscr{F})$ is a
\emph{Gramian} matrix and hence positive semidefinite. Therefore
adding the small regularization term $\lambda^{(t)} \mathbf{I}$
would make the matrix at the left side always positive definite
and hence the system of equations resulting the updates for
$\boldsymbol{\mu}$ always has a unique solution. The left hand
side matrix being strictly positive definite guarantees
$\Delta\boldsymbol{\mu}^{(t)}=\boldsymbol{\mu}^{(t+1)}-\boldsymbol{\mu}^{(t)}$
to be in the descent direction since we have
\begin{equation}
\label{eq31d} [\Delta\boldsymbol{\mu}^{(t)}]^{\mathbf{T}}
[\mathbf{J}_{\boldsymbol{\mu}}(\mathscr{F})]=-[\Delta\boldsymbol{\mu}^{(t)}]^{\mathbf{T}}
[ \mathbf{\tilde{H}}_{\boldsymbol{\mu}}(\mathscr{F})+\lambda^{(t)}
\mathbf{I}][\Delta\boldsymbol{\mu}^{(t)}]<0
\end{equation}
where the superscript ${}^\mathbf{T}$ denotes the matrix
transpose. More technical details about the implementation of the
Levenberg-Marquardt algorithm such as the techniques of choosing
$\lambda^{(t)}$ at each iteration based on trust region algorithms
are available in
\cite{dennis1996numerical,bertsekas1999nonlinear}.

We now turn our attention to the determination of
$p_i(\mathbf{x})$ and $p_o(\mathbf{x})$.  For simplicity here we
assume that $p(\mathbf{x})$ is piecewise constant as is often the
case in work of this type \cite{tai2004survey}. This, in addition
to the shape parameters, we need to determine two constants
defining the contrasts over the regions. We do note that the
approach can easily be extended to consider other low order
``texture models'' as is done in e.g. \cite{kilmer20033d}. As our
primary interest in this paper is a new method for representing
shape, we defer this work to the future. Under our piecewise
constant contrast model, we denote the contrasts as
$p_i(\mathbf{x})=\mathpzc{p}_i$ and
$p_o(\mathbf{x})=\mathpzc{p}_o$, were $\mathpzc{p}_i$ and
$\mathpzc{p}_o$ are unknown constant values. Following analogous
equations as (\ref{eq28}) and (\ref{eq30}), these parameters can
also be appended to the unknown PaLS parameters. The sensitivity
of $\mathscr{F}$ with respect to these parameters can also be
derived based on the fact that
\begin{equation}
\label{eq32} \frac{\partial p}{\partial
\mathpzc{p}_i}=1-\frac{\partial p}{\partial
\mathpzc{p}_o}=H_{rg}(\phi-c).
\end{equation}
and the second order derivatives of $p$ with respect to
$\mathpzc{p}_i$ and $\mathpzc{p}_o$ are zero.

For the PaLS approach represented in this section, the intention
is to keep the formulations general and emphasize the fact that
this representation can formally reduce the dimensionality of the
shaped based inversion. Clearly the expression $\partial
\phi/\partial \mu_j$ depends on how the PaLS functions are related
to their parameters. In the next section a specific PaLS
representation is presented the efficiency of which will be later
examined through some examples.

\section{PaLS Function Representation}\label{sec5}
\subsection{Adaptive Radial Basis Functions}
As pointed out earlier, appropriate choice of a PaLS function can
significantly reduce the dimensionality of an inverse problem. In
this paper we are interested in a low order model for the level
set function which will provide flexibility in terms of its
ability to represent shapes of varying degree of complexity as
measured specifically by the curvature of the boundary. This
representation may allow for ``coarse scale'' elements capable of
representing boundaries of low curvature with few elements.
Furthermore, it is desired to have finer grain elements capable of
capturing higher curvature portions of the boundary such as sharp
turns or corners. Such adaptability is desirable for a general
PaLS representation since the characteristics must be present for
well-posed inverse problems such as full view X-ray CT or even
image segmentation where high fidelity reconstructions are
possible. On the other hand, as we demonstrate in Section
\ref{sec6}, for severely ill-posed problems, the availability of
models with this parsimonious, but flexible structure may allow
for the recovery of geometric detail that otherwise would be not
be obtainable form e.g., a traditional level set or pixel based
approach. Writing a PaLS function as a weighted summation of some
basis functions may be a reasonable choice here, where different
terms in the summation may handle some desired properties about
the shape. Here we focus specifically on the class of radial basis
functions (RBF). We are motivated to concentrate on RBFs as they
have shown to be very flexible in representing functions of
various detail levels. This flexibility makes them appropriate
choice for topology optimization \cite{wang2006radial}, solving
partial differential equations \cite{kansa1990multiquadrics} and
multivariate interpolation of scattered data
\cite{hardy1971multiquadric,wendland2005scattered}. Some examples
of commonly used RBFs are Gaussian, multiquadric, polyharmonic
splines and thin plate splines. More details about these functions
and their properties are available in \cite{buhmann2001radial}.

Based on the statements made, consider the PaLS function
\begin{equation}
\label{eq33}
\phi(\mathbf{x},\boldsymbol{\alpha})=\sum_{j=1}^{m_0}\alpha_j\psi(\|\mathbf{x}-\boldsymbol{\chi}_j
\|)
\end{equation}
where $\psi:\mathbb{R}^{+}\rightarrow \mathbb{R}$ is a
sufficiently smooth radial basis function, $\boldsymbol{\alpha}
\in \mathbb{R}^{m_0}$ is the PaLS parameter vector and $\|.\|$
denotes the Euclidean norm. The points $\boldsymbol{\chi}_j$ are
called the RBF centers. In an interpolation context, usually the
centers are decided in advance and distributed more densely in
regions with more data fluctuations. However, in a PaLS approach
there may be limited information about the shape geometry and
therefore more flexibility is required. Thus here we consider a
more general PaLS function of the form
\begin{equation}
\label{eq34}
\phi(\mathbf{x},[\boldsymbol{\alpha},\boldsymbol{\beta},\boldsymbol{\chi}])=\sum_{j=1}^{m_0}\alpha_j\psi(\|\beta_j(\mathbf{x}-\boldsymbol{\chi}_j)
\|^\dag),
\end{equation}
for which the vector of centers
$\boldsymbol{\chi}=[\boldsymbol{\chi}_1,\boldsymbol{\chi}_2,\cdots,\boldsymbol{\chi}_{m_0}]$
and the dilation factors
$\boldsymbol{\beta}=[\beta_1,\beta_2,\cdots,\beta_{m_0}]$ are
added to the set of PaLS parameters. Also in order to make the
PaLS function globally differentiable with respect to the elements
of $\boldsymbol{\beta}$ and $\boldsymbol{\chi}$, similar to
\cite{acar1994analysis} a smooth approximation of the Euclidean
norm is used as
\begin{equation}
\label{eq35} \|\mathbf{x} \|^\dag:=\sqrt{\|\mathbf{x}
\|^2+\upsilon^2} \qquad \forall \mathbf{x} \in \mathbb{R}^n,
\end{equation}
where $\upsilon\neq 0$ is a small real number. The use of
(\ref{eq34}) rather than (\ref{eq33}) makes the PaLS function
capable of following more details through scaling the RBFs or
floating centers moving to regions where more details are
required. To incorporate this into the optimization methods
described, the sensitivities of $\phi$ with respect to the PaLS
parameters are
\begin{equation}
\label{eq36} \frac{\partial \phi}{\partial
\alpha_j}=\psi(\|\beta_j(\mathbf{x}-\boldsymbol{\chi}_j) \|^\dag)
\end{equation}
and
\begin{equation}
\label{eq37} \frac{\partial \phi}{\partial
\beta_j}=\alpha_j\beta_j\frac{\|(\mathbf{x}-\boldsymbol{\chi}_j)
\|^2}{\|\beta_j(\mathbf{x}-\boldsymbol{\chi}_j)
\|^\dag}\psi'(\|\beta_j(\mathbf{x}-\boldsymbol{\chi}_j) \|^\dag).
\end{equation}
Also considering $\chi_j^{(k)}$ and $x^{(k)}$ to be the
$k^{\mbox{th}}$ components of $\boldsymbol{\chi}_j$ and
$\mathbf{x}$ as points in $\mathbb{R}^n$, for $k=1,2,\cdots,n$ we
have
\begin{equation}
\label{eq38} \frac{\partial \phi}{\partial
\chi_j^{(k)}}=\alpha_j\beta_j^2\frac{\chi_j^{(k)}-x^{(k)}}{\|\beta_j(\mathbf{x}-\boldsymbol{\chi}_j)
\|^\dag}\psi'(\|\beta_j(\mathbf{x}-\boldsymbol{\chi}_j) \|^\dag).
\end{equation}
Clearly the sensitivities obtained are general and valid for any
RBF in $C^1(\mathbb{R}^+)$. In the next section, we consider a
specific class of RBFs that we have found to be particularly well
suited to the PaLS problem.

\subsection{The Choice of Compactly Supported Radial Basis Functions}
Usually the RBFs used in various applications involving function
representations are $C^\infty$ functions with global support
\cite{buhmann2001radial}. However a different class of RBFs which
have recently been under consideration are the compactly supported
radial basis functions (CSRBFs) \cite{wendland2005scattered}.
These functions become exactly zero after a certain radius while
still retaining various orders of smoothness. From a numerical
point of view, compact support of the RBFs yields sparsity in the
resulting matrices arising in the implementation of these methods
and hence reduces the computation cost. This was recently the
motivation to use these functions in simplifying level set methods
\cite{gelas2007compactly}. Another interesting property of these
functions is their local sensitivities to the underlying
parameters \cite{wendland2005scattered}. In other words, when a
function is expressed as a weighted sum of CSRBFs, changing a term
would not have a global effect on the function and only locally
deforms it.

\begin{figure}[t]
\hspace{-.75cm}
\begin{tabular}{ccc}
\epsfig{file=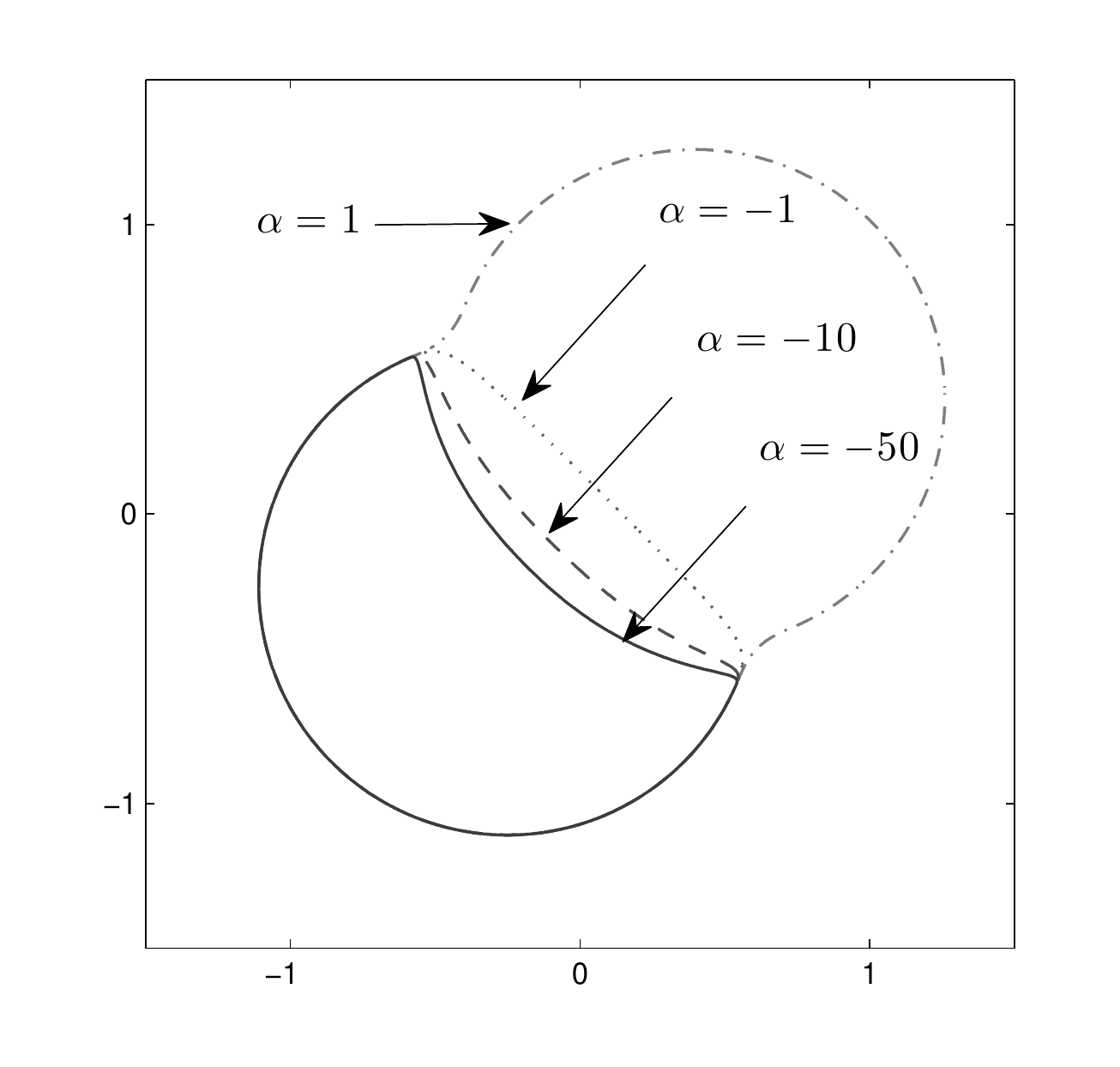,width=0.35\linewidth,clip=a}
\epsfig{file=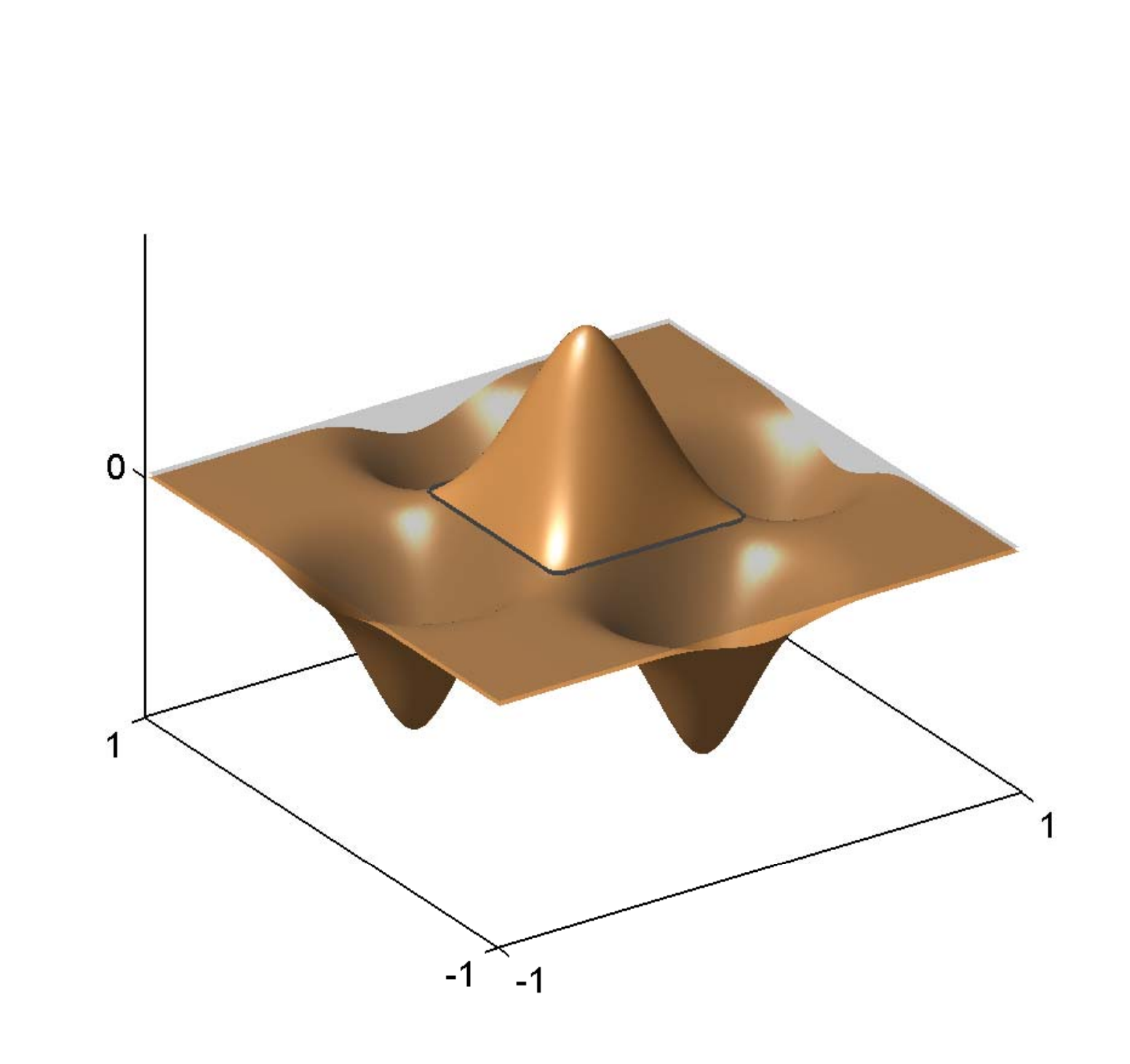,width=0.35\linewidth,clip=a}
\epsfig{file=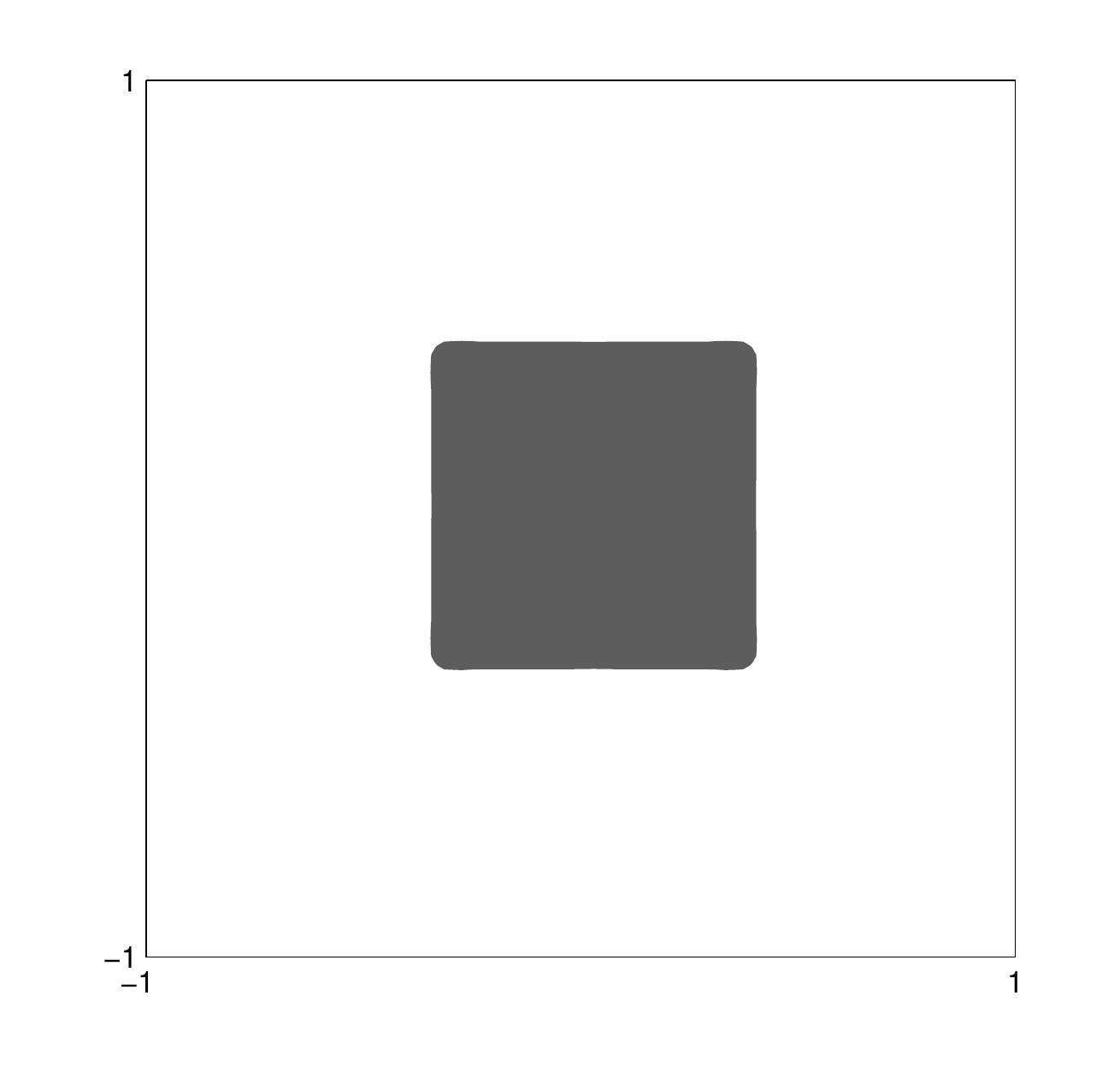,width=0.35\linewidth,clip=}
\end{tabular}
\caption{(a) Left: A close to zero ($c=0.01$) level set of the
function $\psi_1+\alpha \psi_2$ for various values of $\alpha$.
The CSRBF used is the Wendland's function $\psi_{1,1}$ (cf. Table
\ref{table2}), the dilation factors are $\beta_1=\beta_2=1$ and
the centers are taken as $\chi_1=(-\frac{1}{4},-\frac{1}{4})$ and
$\chi_2=(\frac{2}{5},\frac{2}{5})$. (b) Center: A PaLS function
involved in representation of a square at a close to zero level
set. Only 5 bumps are involved. (c) Right: The representation of
the square}\label{fig2}
\end{figure}

Beside aforementioned advantages of the CSRBFs, our interest in
this class of RBFs arises from their potential in reconstructing
the shapes with a very small number of terms in a PaLS
representation as (\ref{eq34}). Furthermore, as will be explained,
this representation can involve shapes with corners and rather
high curvature regions. For $\psi\geq 0$ being a smooth CSRBF,
lets denote every basis term in (\ref{eq34}) as
\begin{equation}
\label{eq40}
\psi_j(\mathbf{x})=\psi(\|\beta_j(\mathbf{x}-\boldsymbol{\chi}_j)
\|^\dag)
\end{equation}
and call $\psi_j$ a \emph{bump}. Due to the compact support of
these functions, for every two bumps $\psi_j$ and $\psi_k$ we can
write
\begin{equation}
\label{eq41}
\supp(\psi_j+\psi_k)=\supp(\psi_j)\;\cup\;\;\supp(\psi_k).
\end{equation}
For a real valued function $\vartheta$ defined over $\mathbb{R}^n$
we define
\begin{equation}
\label{eq41.a} \mathscr{I}\!\!_c(\vartheta):=\{\mathbf{x} :
\vartheta(\mathbf{x})\geq c \}.
\end{equation}
Clearly for $c> 0$, $\mathscr{I}_c(\vartheta)$ represents the
interior of the $c$-level set of $\vartheta$. Based on
(\ref{eq41}) and the smoothness of the bumps, we obviously have
that as $c\rightarrow 0^+$, $\mathscr{I}\!\!_c(\psi_j+\psi_k)$
would tend to
$\mathscr{I}\!\!_c(\psi_j)\bigcup\mathscr{I}\!\!_c(\psi_k)$. More
generally for $\alpha_j>0$,
$\mathscr{I}\!\!_c(\sum_{j=1}^{m_0}\alpha_j\psi_j)$ tends to
$\bigcup_{j=1}^{m_0} \mathscr{I}\!\!_c(\psi_j)$ as
$(c/\alpha_j)\rightarrow 0^+$. In other words, using relatively
large CSRBF weights (as compared to $c$) can imply reconstruction
of the shape through the union of a collection of floating balls
of various radii. Moreover $\mathscr{I}\!\!_c(\psi_j-\alpha_k
\psi_k)$ would tend to $\mathscr{I}\!\!_c(\psi_j)\! \setminus\!
\mathscr{I}\!\!_c(\psi_k)$ as $c\rightarrow 0^+$ and $\alpha_k
\rightarrow +\infty$. In a more general fashion,
$\mathscr{I}\!\!_c(\alpha_j \psi_j-\alpha_k \psi_k)$ tends to
$\mathscr{I}\!\!_c(\psi_j)\! \setminus\!
\mathscr{I}\!\!_c(\psi_k)$ as $(c/\alpha_j)\rightarrow 0^+$ and
$(\alpha_k/\alpha_j) \rightarrow +\infty$. Therefore in this
context, bumps with larger negative coefficients can yield holes
or inflect the shape by excluding some portions of it. We would
consider the two aforementioned properties of the CSRBFs as a
``pseudo-logical'' behavior of these functions. This property can
result in rather high curvature geometries with a limited number
of bumps. Besides high curvature regions, low curvature segments
(e.g., an almost straight line in $\mathbb{R}^2$ or a planar
segment in $\mathbb{R}^3$) can be formed by interaction of two
identical bumps but with opposite signs at their footprint
intersections. Figure \ref{fig2}.a sheds more light on
aforementioned facts, and shows the interaction of two bumps,
which for instance can represent a crescent with two rather sharp
corners ($\alpha=-50$), or a contour with a low curvature segment
($\alpha=-1$). As a second example, a representation of a square
using only 5 bumps is depicted in Figures \ref{fig2}.b and
\ref{fig2}.c.

The most commonly used CSRBFs are those called Wendland's
functions \cite{wendland2005scattered}. The smoothness and the
compact support provided by Wendland's functions are of interest
and hence we shall use them as the basis for our PaLS approach.
Wendland's functions follow the general form of
\begin{equation}
\label{eq39} \psi_{n,l}(r)=\left\{
     \begin{array}{lr}
       P_{n,l}(r) & 0\leq r \leq 1\\
       0 & r>1
     \end{array}
   \right.
\end{equation}
when representing an RBF in $\mathbb{R}^n$, with $P_{n,l}$ being a
univariate polynomial of degree $\lfloor n/2\rfloor +3l+1$. In
terms of smoothness, this class of RBFs belong to $C^{2l}$. A
derivation of these functions is provided in
\cite{wendland2005scattered}, from which we have listed the first
few functions in Table \ref{table2}.
\begin{table}[t]
\caption{Compactly supported RBFs of minimal degree where $\ell=\lfloor n/2 \rfloor +l+1$ and $(1-r)_+=\max\{0,1-r\}$} 
\centering 
\begin{tabular}{c| c} 
\hline\hline
Function & Smoothness \\ [0.5ex] 
\hline 
$\psi_{n,1}(r)=(1-r)_+^{\ell+1}((\ell+1)r+1)$ & $C^2$\\
$\psi_{n,2}(r)=(1-r)_+^{\ell+2}((\ell^2+4\ell+3)r^2+(3\ell+6)r+3)$ & $C^4$\\
$\psi_{n,3}(r)=(1-r)_+^{\ell+3}((\ell^3+9\ell^2+23\ell+15)r^3+$ & $C^6$\\
$\hspace{2cm}(6\ell^2+36\ell+45)r^2+(15\ell+45)r+15)$ & \\
 \hline
\end{tabular}
\label{table2} 
\end{table}

In the next section we discuss the regularized heaviside function
and explain how choosing an appropriate version of this function
in (\ref{eq25}), can pave the path for exploiting the
pseudo-logical behavior of the bumps.

\subsection{Numerical Approximation of the Heaviside Function}

In a shape-based representation such as (\ref{eq25}), solving the
inverse problem numerically and making the evolution of the level
set function possible requires using a smooth version of the
heaviside function. A possible $C^\infty$ regularization of $H(.)$
denoted as $H_{1,\epsilon}$ is the one used in
\cite{chan2001active}, as
\begin{equation}
\label{eq42} H_{1,\epsilon}(x)=\frac{1}{2}\big (
1+\frac{2}{\pi}\arctan(\frac{\pi x}{\epsilon})  \big).
\end{equation}
This function is commonly used in shape-based applications and
specifically in the recent parametric representations of the level
set function in image segmentation \cite{gelas2007compactly,
bernard2009variational}. Chan {\it et al.} in
\cite{chan2006algorithms} have studied the characteristics of the
resulting level set functions using $H_{1,\epsilon}$ as the
evolution proceeds. In the context of the current shape-based
problem for which $p_i(\mathbf{x})=\mathpzc{p}_i$ and
$p_o(\mathbf{x})=\mathpzc{p}_o$, and considering the zero level
set, (\ref{eq7}) reveals that in an evolution process $\phi$ is
likely to evolve towards a state that $H(\phi(\mathbf{x}))= 1$ for
$\mathbf{x} \in D$ and $H(\phi(\mathbf{x}))= 0$ for $\mathbf{x}
\in \Omega \setminus D$.
\begin{figure}[t]
\hspace{-.55cm}
\begin{tabular}{cc}
\epsfig{file=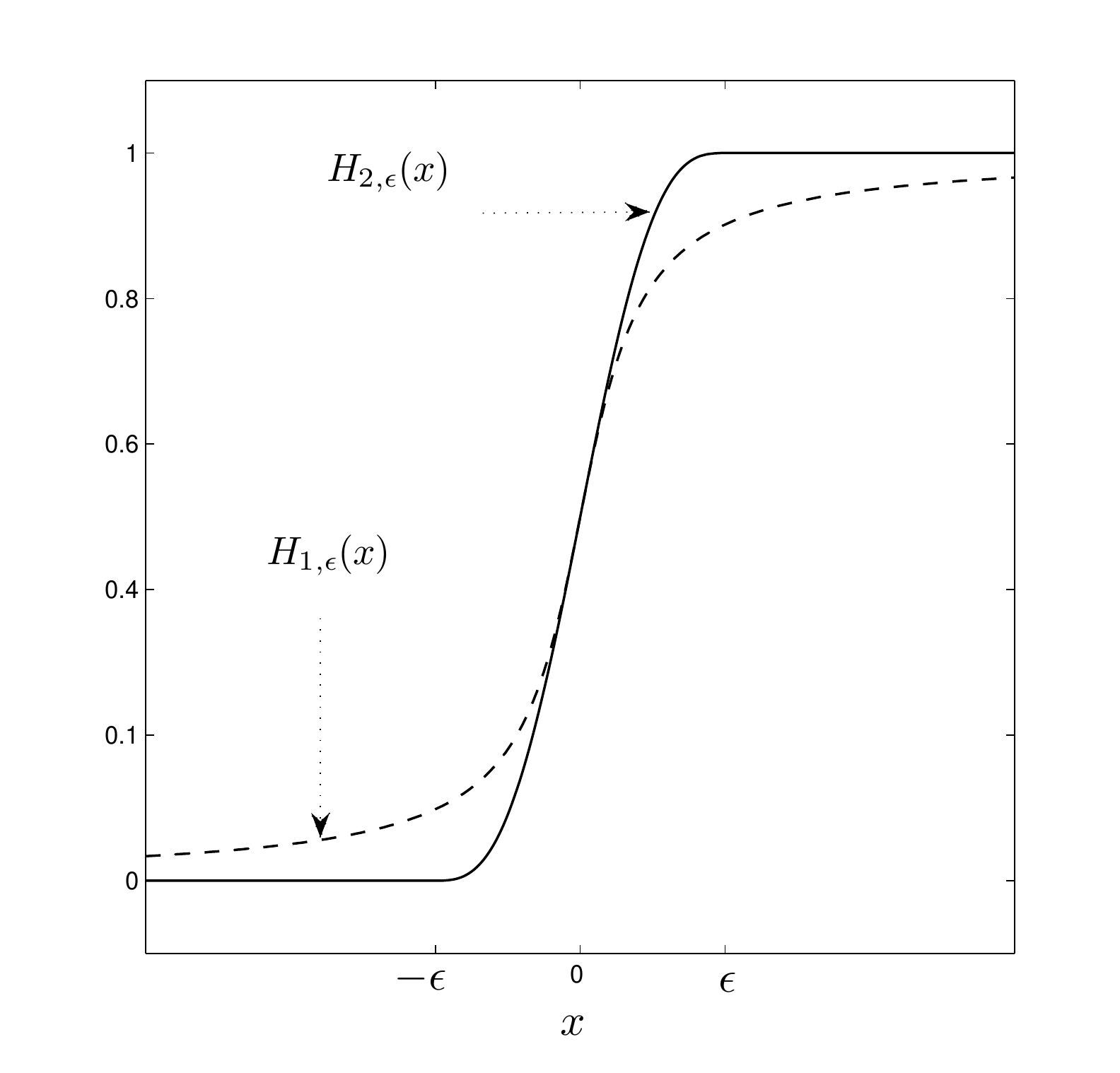,width=0.5\linewidth,clip=a}&
\epsfig{file=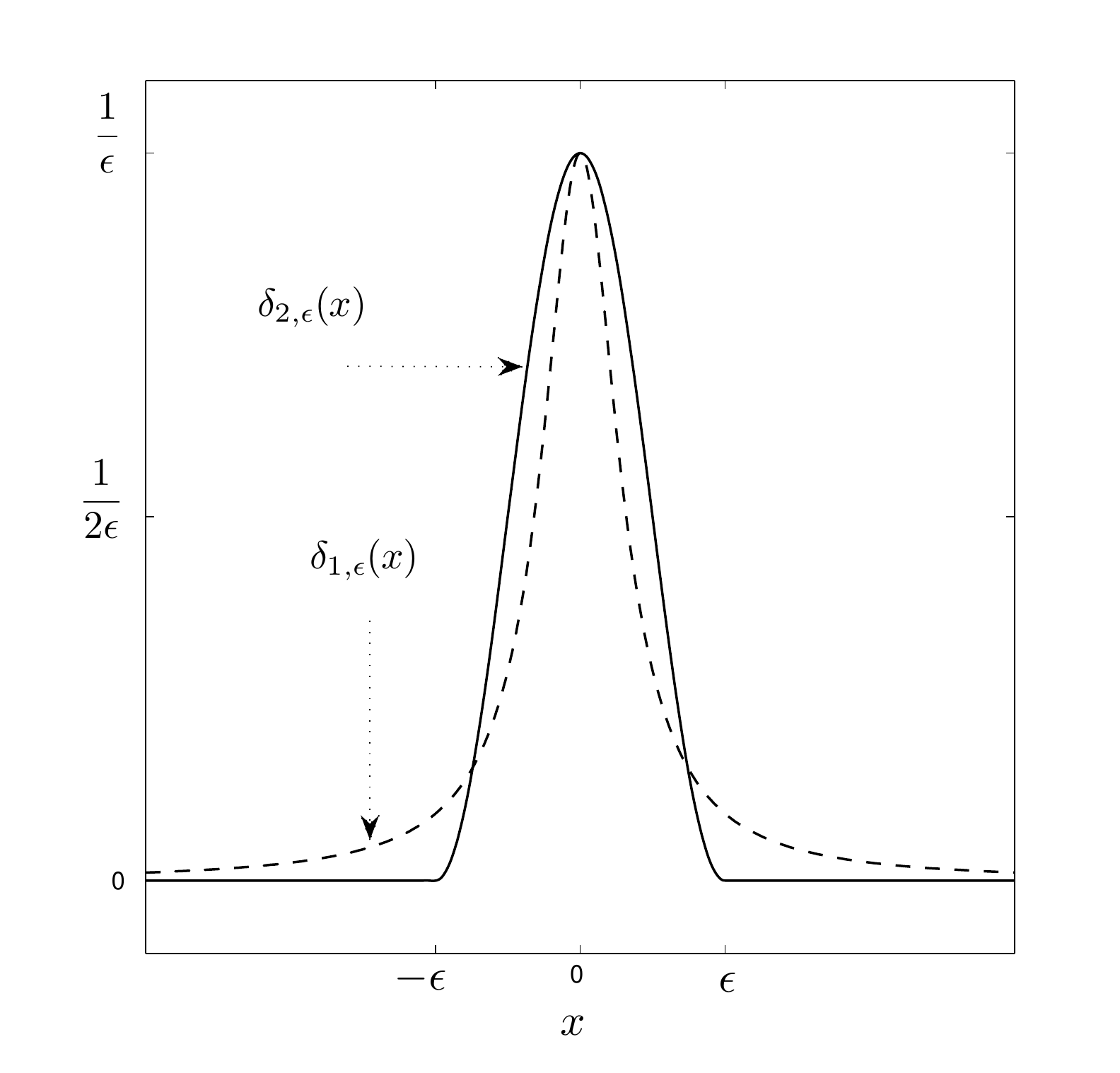,width=0.5\linewidth,clip=}
\end{tabular}
\caption{(a) Left: Two regularized versions of the heaviside
function. (b) Right: Corresponding regularized delta
functions}\label{fig3}
\end{figure}
Referring to Figure \ref{fig3}.a, one observes that when
$H_{1,\epsilon}$ is used as the regularized heaviside function, in
order to have $H_{1,\epsilon}(\phi)\simeq 1$ in $D$, the level set
function should take rather large positive values (as compared to
$\epsilon$) in this region. Analogously in $\Omega \setminus D$,
$\phi$ is pushed to take rather large negative values. These
constraints are implicitly imposed on the resulting level set
function and specifically using CSRBFs in a PaLS approach as
(\ref{eq34}), the bumps are expected to distribute throughout
$\Omega$ to form a level set function which takes rather large
positive (or negative) values inside (or outside) $D$.

\begin{figure}[t]
\hspace{-.85cm}
\begin{tabular}{cc}
\epsfig{file=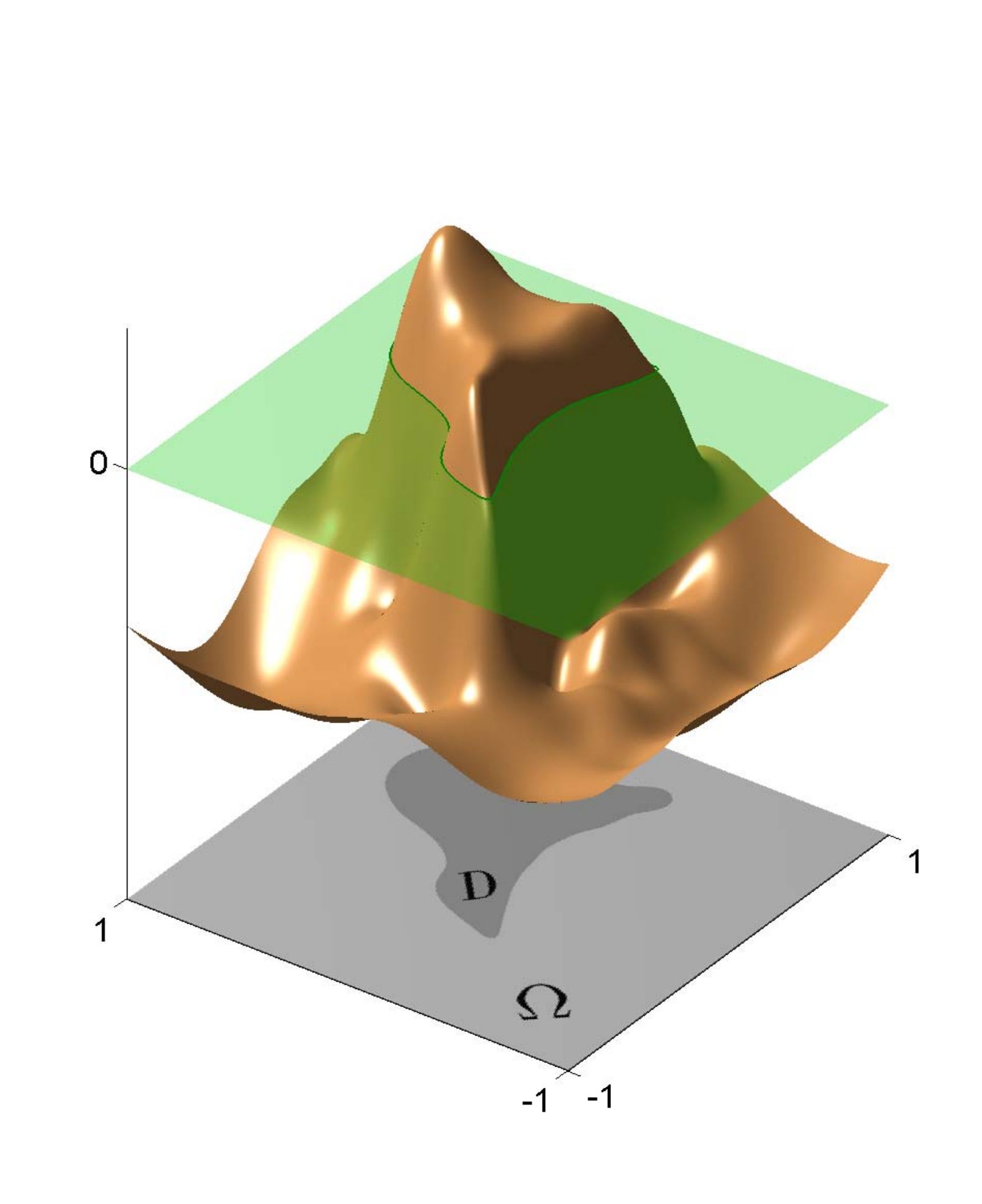,width=0.52\linewidth,clip=a}&
\epsfig{file=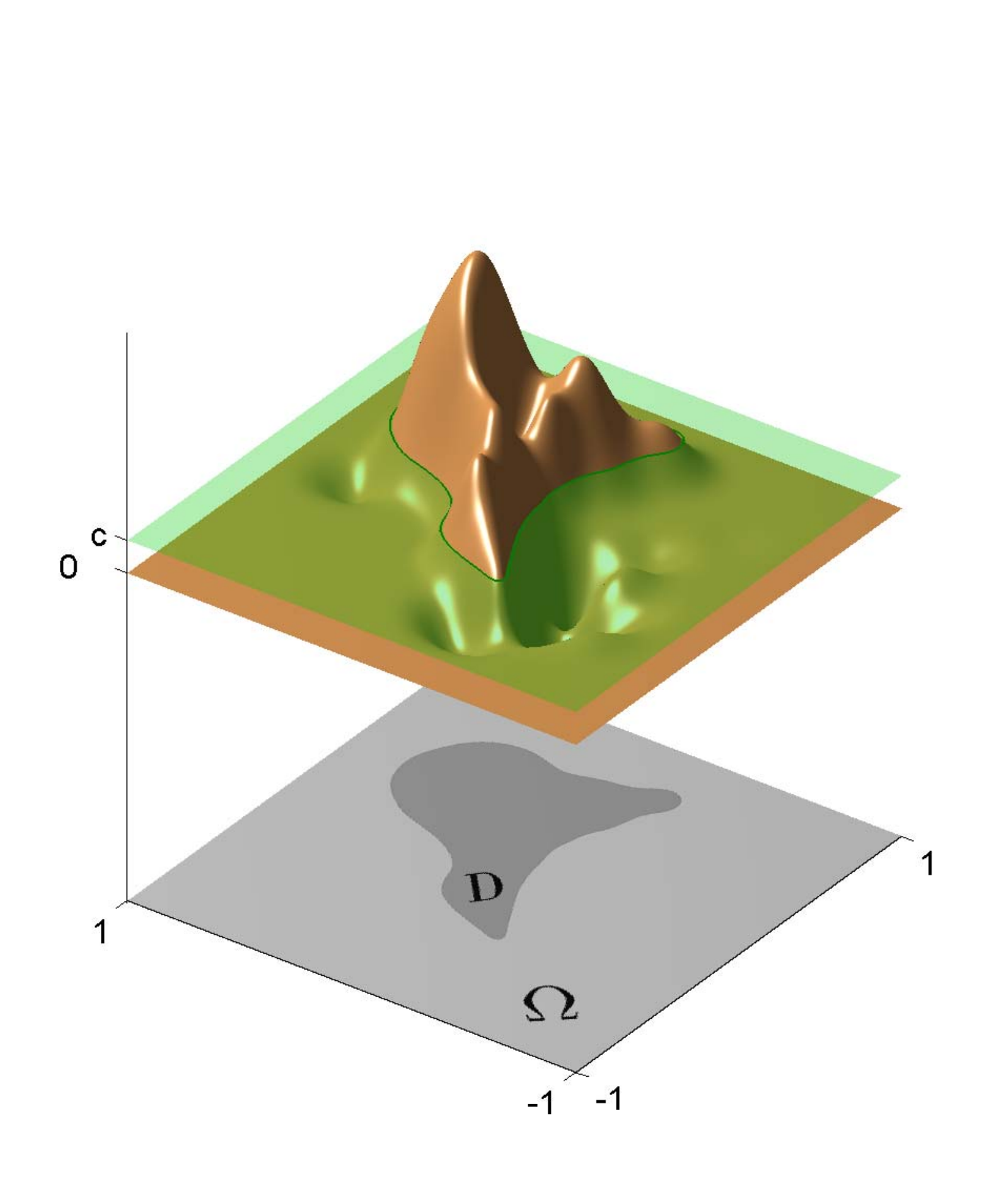,width=0.52\linewidth,clip=}
\end{tabular}
\caption{(a) Left: A typical PaLS function resulted using
$H_{1,\epsilon}$ with 65 bumps and considering the zero level set.
(b) Right: A typical PaLS function resulted using $H_{2,\epsilon}$
with 30 bumps and considering the $c$ level set}\label{fig4}
\end{figure}

An alternative choice of the regularized Heaviside function is the
$C^2$ function,
\begin{equation}
\label{eq43} H_{2,\epsilon}(x)= \left \{
\begin{array}{cl}
\hspace{-2cm}1 &  x>\epsilon \\
\hspace{-2cm} 0 &  x<-\epsilon \\
\frac{1}{2}+\frac{x}{2\epsilon}+\frac{1}{2\pi}\sin(\frac{\pi x}{\epsilon}) &  |x|\leq \epsilon,\\
\end{array} \right.
\end{equation}
as proposed in \cite{zhao1996variational}. It can be shown that if
the nonzero level set $c$ and $\epsilon$ are appropriately set,
using $H_{2,\epsilon}$ enables us to exploit the pseudo-logical
behavior described previously. More specifically, for $c>0$, and
$\phi$ being a weighted sum of some bumps, and therefore compactly
supported itself, we clearly want the points for which $\phi\leq
0$ to belong to $\Omega \setminus D$, i.e., to correspond to the
intensity $\mathpzc{p}_o$. Using (\ref{eq25}), this requires
having $H(\phi-c)=0$ for $\phi\leq 0$. Referring to Figure
\ref{fig3}.a, this condition is satisfied when $H_{2,\epsilon}$ is
chosen as the regularized heaviside function and $-c\leq
-\epsilon$ (or in general case $|c|\geq \epsilon$ as the required
criteria). Practically, this use of $H_{2,\epsilon}$ takes away
the implicit constraint imposed on the level set function using
$H_{1,\epsilon}$, i.e., the bumps do not have to spread throughout
$\Omega$ and may only concentrate inside and about the shape $D$
to perform a higher resolution shaping. Figure \ref{fig4} shows a
typical shape representation resulted using $H_{1,\epsilon}$ and
$H_{2,\epsilon}$ highlighting the pseudo-logical property.

\begin{figure}[t]
\hspace{-.85cm}
\begin{tabular}{cc}
\epsfig{file=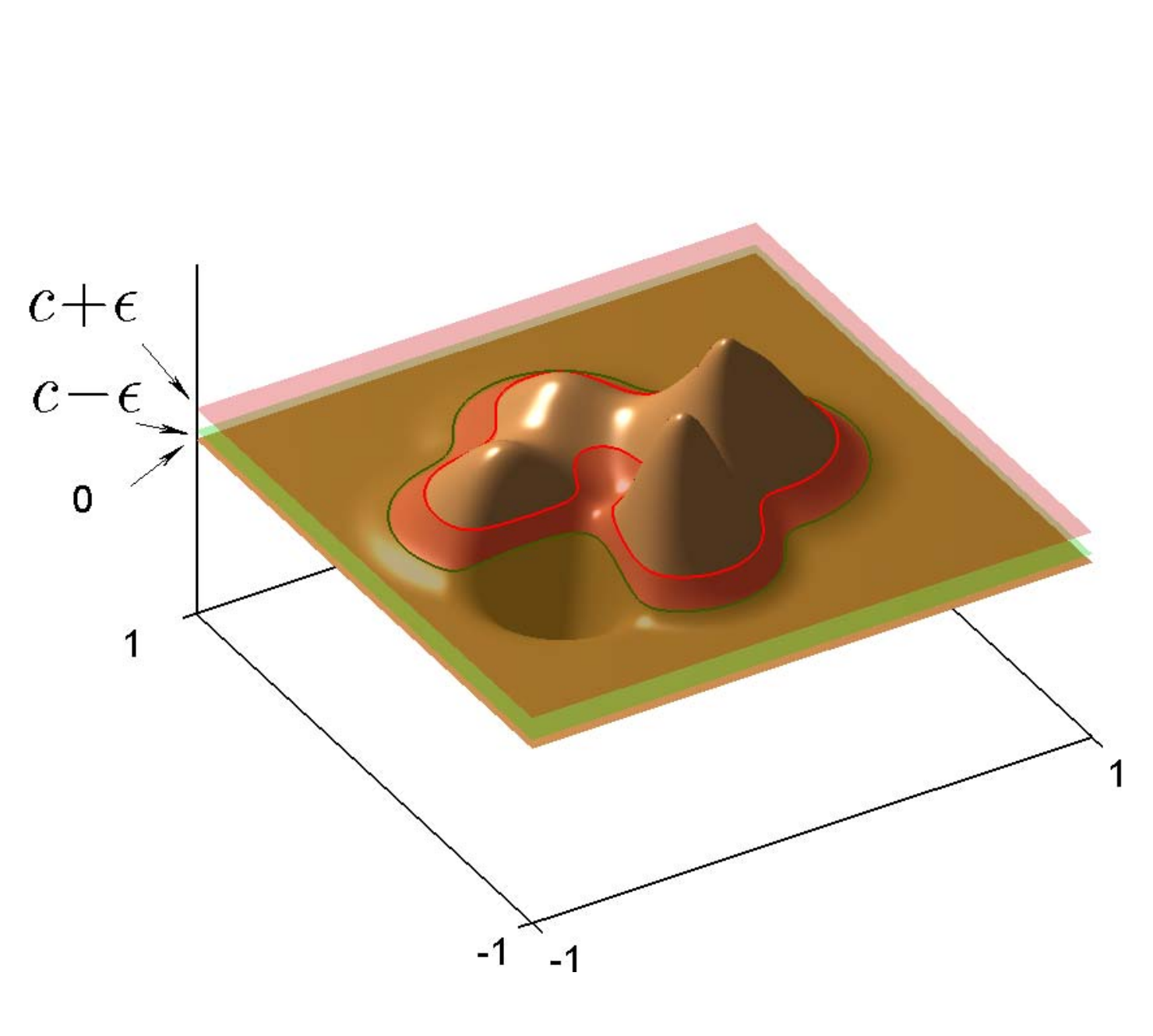,width=0.52\linewidth,clip=a}&
\epsfig{file=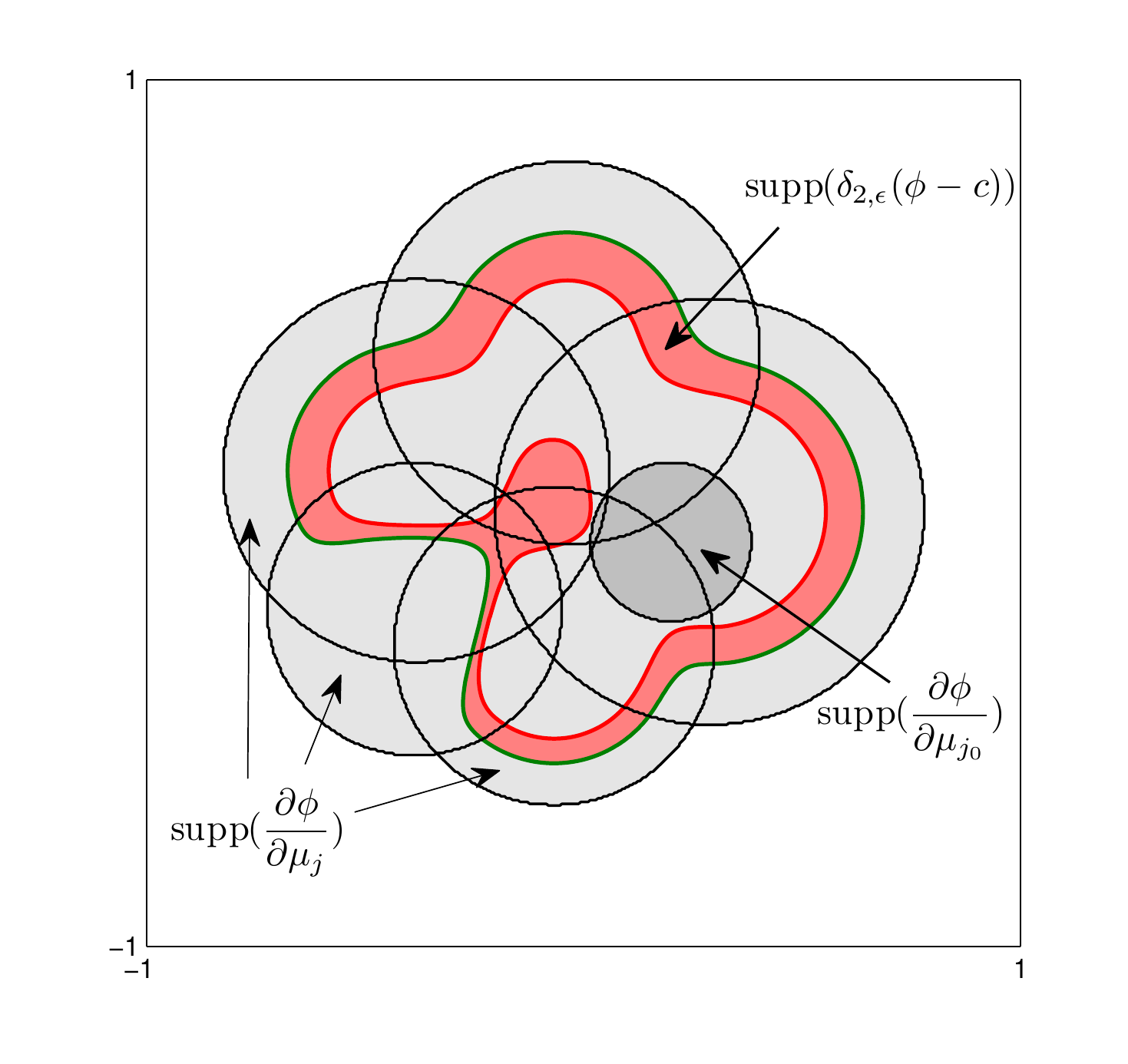,width=0.52\linewidth,clip=}
\end{tabular}
\caption{(a) Left: A typical PaLS function composed of 6 bumps,
and the $c\pm \epsilon$ level sets.\hspace{3cm} (b) Right: The
bumps with active evolving parameters and the frozen
bump}\label{fig5}
\end{figure}

Besides exploiting the pseudo-logical behavior, using
$H_{2,\epsilon}$ can provide another advantage which reduces the
dimensionality of the problem at every iteration, i.e., at every
step of the evolution process, the cost function will be sensitive
to a specific group of PaLS parameters. To further describe this
behavior, consider the term $\delta_{rg}(\phi-c)\frac{\partial
\phi } {\partial \mu_j}$ appearing in both the Jacobian and the
approximate Hessian of $\mathscr{F}$ in (\ref{eq28.x}) and
(\ref{eq31.b}). Also for an evolution iteration, assume the
superscript $(t)$ indicating the value of every quantity at that
iteration. If for a PaLS parameter such as $\mu\!_{j_{{}_0}}$
\begin{equation}
\label{eq44} \delta_{rg}(\phi^{(t)}-c)\frac{\partial \phi^{(t)} }
{\partial \mu\!_{j_{{}_0}}}\Big
|_{\mu\!_{j_{{}_0}}=\mu\!_{j_{{}_0}}^{(t)}} = 0 \qquad \forall
\mathbf{x}\in \Omega,
\end{equation}
then using either one of the minimization schemes as (\ref{eq28a})
or (\ref{eq31c}) yields
\begin{equation}
\label{eq45}  \mu\!_{j_{{}_0}}^{(t+1)}= \mu\!_{j_{{}_0}}^{(t)}.
\end{equation}
The reason for this result is clear for the gradient descent
scheme (\ref{eq28a}), based on the fact that using (\ref{eq44}) in
(\ref{eq28.x}) causes the $j_{{}_0}^{th}$ element of gradient
vector to vanish and hence $\mu\!_{j_{{}_0}}^{(t)}$ remaining
unchanged at the corresponding iteration. For the
Levenberg-Marquardt scheme (\ref{eq31c}), using (\ref{eq44}) in
(\ref{eq31.b}) causes all the elements in the $j_{{}_0}^{th}$ row
of the approximate Hessian matrix to vanish, which results the
corresponding equation
\begin{equation}
\label{eq46} \lambda(\mu\!_{j_{{}_0}}^{(t+1)}-
\mu\!_{j_{{}_0}}^{(t)}) =0,
\end{equation}
again equivalent to (\ref{eq45}). Therefore, for either one of the
minimization schemes, if (\ref{eq44}) holds, the PaLS parameter
$\mu\!_{j_{{}_0}}$ will stay unchanged in that iteration. We here
describe a common case that (\ref{eq44}) holds during the
evolution process:

By using $H_{2,\epsilon}$ as the regularized level set function,
the corresponding regularized delta function $\delta_{2,\epsilon}$
will be compactly supported (as shown in Figure \ref{fig3}.b),
hence $\delta_{2,\epsilon}(\phi-c)$ is only nonzero for
$c-\epsilon<\phi<c+\epsilon$. On the other hand, based on the PaLS
approach presented in this paper using CSRBFs, for $\mu_j$ being
any of the PaLS parameters $\alpha_j$, $\beta_j$ or $\chi_j^{(k)}$
corresponding to the bump $\psi_j$, we have
\begin{equation}
\label{eq47} \supp(\frac{\partial \phi } {\partial
\mu_j})\subseteq \supp(\psi_j).
\end{equation}
This fact is easily observable in (\ref{eq36}), (\ref{eq37}) and
(\ref{eq38}), where the related derivatives can only have nonzero
values in $\supp(\psi_j)$. Therefore, if at some iteration and for
a bump $\psi_{j_{{}_0}}$,
\begin{equation}
\label{eq48} \supp\big(\delta_{2,\epsilon}(\phi-c)\big) \cap
\supp(\psi_{j_{{}_0}})=\emptyset,
\end{equation}
then in that iteration we have
\begin{equation}
\label{eq49} \delta_{2,\epsilon}(\phi-c)\frac{\partial \phi }
{\partial \mu\!_{j_{{}_0}}}=0
\end{equation}
and therefore the PaLS parameters corresponding to
$\psi_{j_{{}_0}}$ will stay unchanged in that iteration. Figure
\ref{fig5} illustrates this phenomenon, showing a PaLS function
composed of 6 bumps at some iterations. For 5 of the bumps used,
the corresponding support does intersect the region
$\supp(\delta_{2,\epsilon}(\phi-c))$, and therefore their
corresponding parameters have the potential to change at this
state of the PaLS function. However, a bump denoted as
$\psi_{j_{{}_0}}$, does not intersect
$\supp\big(\delta_{2,\epsilon}(\phi-c))$, and the underlying PaLS
parameters do not need to be considered in that iteration. This
approach is similar to the \emph{narrow-banding} approach in
traditional level-set methods
\cite{adalsteinsson1995fast,sethian1999level,peng1999pde}, where
the values of the level set function are only updated on a narrow
band around the zero level set and hence reducing the computation
load. In our approach, however, this band is the points for which
$c-\epsilon<\phi<c+\epsilon$ and the bumps which do not intersect
with this band do not evolve at the corresponding iteration and
hence their corresponding parameters are not updated.

In the next section, through a number of examples drawn from a
wide range of applications, we will show the superior performance
of the proposed method specifically exploiting the pseudo-logical
behavior of the CSRBFs.

\section{Examples}\label{sec6}
In this section, we examine our proposed method for three
different inverse problems, namely electrical resistance
tomography, X-ray computed tomography and diffuse optical
tomography. The examples are simulated for 2D imaging and the
results are provided in each section. Throughout all the examples,
$\Omega$ denotes the region to be imaged and $D$ denotes the shape
domain as stated in the previous sections.
\subsection{Electrical Resistance Tomography}

As the first example we consider electrical resistance tomography
(ERT) categorized as a severely ill-posed problem
\cite{hansen1998rank,di2003examples}. The objective of this
problem is the reconstruction of the electrical conductivity
within a region of interest based on the potential or current
measurements performed at the periphery of that region. Such
reconstructions may be applicable in various areas such as medical
imaging \cite{cheney1999electrical}, geophysics
\cite{dines1981analysis} and environmental monitoring
\cite{daily2000electrical}.

For many geophysical applications the underlying physical model
describing the DC potential $u(\mathbf{x})$ inside $\Omega$ in
terms of the conductivity $\sigma(\mathbf{x})$ and the current
source distribution $s(\mathbf{x})$ is
\begin{align}
\label{eq1-1} &\nabla \cdot \big(\sigma(\mathbf{x}) \nabla
u(\mathbf{x})\big)=s(\mathbf{x}) \quad\;\; \mbox{in}\;\; \Omega,\\
\nonumber & \sigma \frac{\partial u}{\partial \nu}=0\qquad
\quad\quad\quad\qquad\quad \mbox{on}\;\;
\partial\Omega_n \subset
\partial \Omega,\\ \nonumber & u=0 \;\;\qquad\quad\quad\quad\qquad\quad\quad \mbox{on}\;\;
\partial\Omega_d=\partial \Omega\setminus\partial\Omega_n,
\end{align}
where $\nu$ denotes the outward unit normal on $\partial \Omega$
and $\partial \Omega_n$ and $\partial \Omega_d$ correspond to
Neumann and Dirichlet boundaries. In many of the applications
(e.g., see
\cite{saunders2005constrained,labrecque1996effects,tripp1984two,
mufti1976finite}) the Dirichlet boundary condition is imposed as
an approximation to the potential in regions far from the actual
imaging region, and is used here for simplicity. For arbitrary
distributions of the conductivity, (\ref{eq1-1}) is usually solved
numerically by means of finite element or finite difference
methods \cite{soleimani2006level,zhang19953}. However, as the main
focus of the paper, we concentrate on piecewise constant
conductivity distribution as $\sigma(\mathbf{x})=\sigma_i$ for
$\mathbf{x} \in D$ and $\sigma(\mathbf{x})=\sigma_o$ for
$\mathbf{x} \in  \Omega\setminus D$.

For the inverse problem, the sensitivities of the measurements to
perturbations of the conductivity (in our approach the
perturbations of the PaLS parameters) are required. For
$s(\mathbf{x})=\delta(\mathbf{x}-\mathbf{x}_s)$, i.e., a point
source current at $\mathbf{x}_s \in \Omega$, we denote by
$u_s(\mathbf{x})$ the resulting potential over the domain $\Omega$
and consider the measured potential at $\mathbf{x}_d \in \Omega$
as
\begin{equation}
\label{eq1-4} u_{ds}=\int_\Omega
u_s(\mathbf{x})\delta(\mathbf{x}-\mathbf{x}_d)\mbox{d}\mathbf{x}.
\end{equation}
The variation of $u_{ds}$ resulting from a perturbation $\delta
\sigma$ in the conductivity (i.e., the Fr\'{e}chet derivative of
the measurements with respect to the conductivity) can be then
expressed as \cite{polydorides2002matlab,
soleimani2008computational}
\begin{equation}
\label{eq1-5} \frac{\mbox{d}u_{ds}}{\mbox{d}\sigma}\;
[\delta\sigma]=\int_\Omega \delta \sigma \;\nabla u_s \cdot \nabla
u_d \; \mbox{d}\mathbf{x},
\end{equation}
where $u_d$ is the adjoint field that results from placing the
current point source at $\mathbf{x}_d$. To express the inverse
problem in a PaLS framework, we consider $\mathcal{M(.)}$ as the
nonlinear forward model mapping the conductivity distribution into
a vector of voltage measurements $\mathbf{u}$ obtained by
performing $M$ experiments, having a different point source
position at each experiment and making $N_\ell$ potential
measurements for $\ell=1,2,\cdots,M$. Having the residual operator
$\mathcal{R}(\sigma)=\mathcal{M}(\sigma)-\mathbf{u}$ and using
(\ref{eq1-5}), the Fr\'{e}chet derivative denoted as
$\mathcal{R}'(\sigma)[.]$ can be considered as a vector consisting
of $M$ sub-vectors $\mathcal{R}'_\ell(\sigma)[.]$, structured as
\begin{equation}
\label{eq1-6} \mathcal{R}'_\ell(\sigma)[\delta\sigma] = \left(
\begin{array}{c}
\int_\Omega \delta \sigma \;\nabla u_\ell \cdot
\nabla u_{1}^{\ell} \; \mbox{d}\mathbf{x}\\
\vdots\\ \int_\Omega \delta \sigma \;\nabla u_\ell \cdot \nabla
u_{N_\ell}^{\ell} \; \mbox{d}\mathbf{x}
\end{array} \right).
\end{equation}
Here $u_\ell$ denotes the potential in the $\ell^{\mbox{th}}$
experiment and $u_i^\ell$ denotes the adjoint field corresponding
to the $\ell^{\mbox{th}}$ experiment resulted from placing the
current point source at the $i^{\mbox{th}}$ measurement point.
Having $\mathcal{R}'(\sigma)[.]$ in hand, one can obtain the PaLS
evolution through using (\ref{eq28}) and (\ref{eq31.b}) in
(\ref{eq31c}).

\begin{figure}[t]
\hspace{-1.5cm}
\begin{tabular}{cc}
\epsfig{file=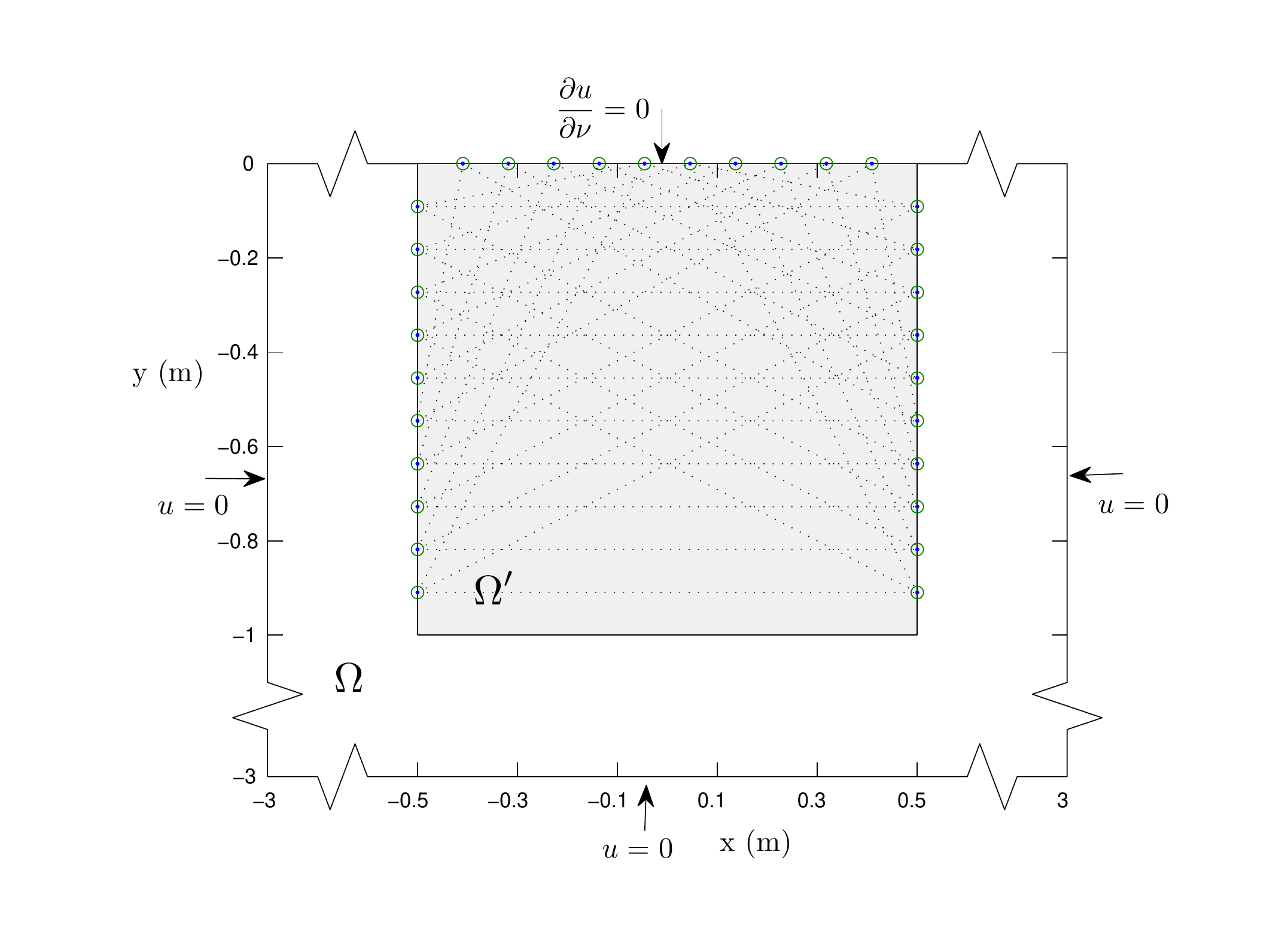,width=0.65\linewidth,clip=a}
\epsfig{file=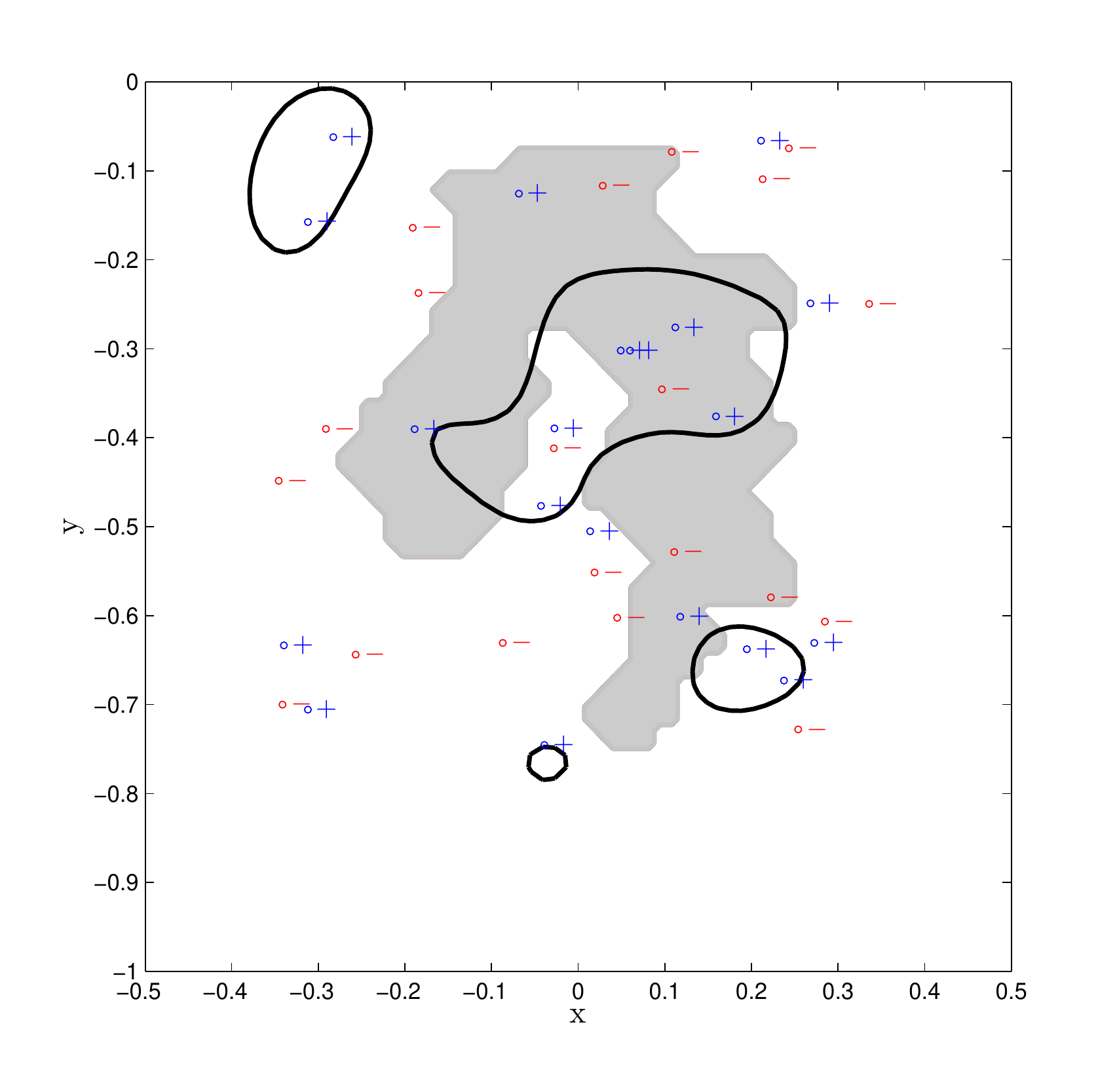,width=0.5\linewidth,clip=a}
\end{tabular}
\caption{(a) Left: The 2D modelling region for ERT. The darker
interior region is the imaging region surrounded by the sources
and detectors. The dashed lines correspond to dipole nodes used in
every experiment (b) Right: The gray region shows the shape to be
reconstructed in the ERT problem. The dots with ``+" and ``-"
signs correspond to the centers of positive and negative weighted
bumps in the initial state of the problem. The black contour is
the resulting $c$-level set of the initial PaLS function.
}\label{figex1}
\end{figure}

For the purpose of this example, we model the electric potential
within the box $\Omega=[-3,3]\times[-3,0]$ all dimensions in units
of meters in $x-y$ plane. Here $\partial\Omega_n$ corresponds to
the top surface ($y=0$) and $\partial\Omega_d$ corresponds to the
sides and bottom boundaries as shown in Figure \ref{figex1}.a. The
region to be imaged is the square
$\Omega'=[-0.5,0.5]\times[-1,0]$, where 30 points are allocated to
the sensors placed equally spaced on the top and the sides of this
region (shown as small circles in the figure). A total of $M=40$
experiments are performed and in every experiment two of the
sensors are used as a current dipole (superposition of two point
sources with opposite signs) and the remaining 28 sensors measure
the potential at the corresponding positions. The dipole sources
are chosen in a cross-medium configuration, where the electrodes
corresponding to every experiment are connected with a dashed line
as again shown in Figure \ref{figex1}.a. With this configuration
we try to enforce electric current flow across the medium
anomalies and obtain data sets more sensitive to shape
characteristics. For the simulation, we use the finite difference
method where $\Omega$ is discretized to 125 grid points in the $x$
direction and 100 points in the $y$ direction. The gridding is
performed in a way that we have a uniform $75\times 75$ grid
within $\Omega'$ (excluding the boundaries containing the sensors)
and the exterior grids linearly get coarser as they get further
from $\Omega'$ in the $x$ and $y$ direction. The forward modelling
is performed over the whole collection of grids in $\Omega$, while
the inversion only involves the pixels within $\Omega'$ (known as
active cells \cite{pidlisecky2007resinvm3d}).

The shape to be reconstructed is shown as the gray region in
Figure \ref{figex1}.b. This shape is a threshholded version of a
real scenario and is of particular interest here because it has a
concavity facing the bottom where there are no measurements
performed. The values for the anomaly conductivities are
$\sigma_i=0.05\;\mbox{Sm}^{-1}$ and
$\sigma_o=0.01\;\mbox{Sm}^{-1}$ and the conductivity value used
for $\Omega\setminus\Omega'$ is the same as the true value of
$\sigma_o$ in all our inversions. In the inversions we consider
the data $\mathbf{u}$ to be the measured potentials generated by
the true anomaly with $1\%$ additive Gaussian noise. For the shape
representation we use $H_{2,\epsilon}$ with $\epsilon=0.1$ and we
consider the $c=0.15$ level set.

For the PaLS representation (\ref{eq33}) we have $m_0=40$ terms
and the Wendland's function $\psi_{1,1}$ is used as the
corresponding bump. For the initial PaLS parameters, we consider a
random initial distribution of the centers $\boldsymbol{\chi}_j$,
within the square $[-0.4,0.4]\times[-0.8,0]$. The weighting
coefficients are initialized as $\alpha_j=\pm 0.2$ where the
centers of alternatively positive and negative weighted bumps are
shown with $``+"$ and $``-"$ signs in Figure \ref{figex1}.b. The
positive initial values of $\alpha_j$ are taken slightly bigger
than $c$ to have some initial $c$-level sets. The purpose behind
having alternative bump signs is to have the narrow-band
$\supp(\delta_{2,\epsilon}(\phi-c))$ cover various regions of
$\Omega'$ and increase the chance of initially involving more
bumps in the shaping as explained in previous section and Figure
\ref{fig5}. The dilation factors are taken uniformly to be
$\beta_j=4$, as an initialization to make the support radii of the
bumps small enough for capturing details and more or less large
enough to carpet the region $\Omega'$. Our intention for this
initialization of the PaLS parameters is to provide a rather
simple, reproducible and general initialization. The stopping
criteria in the reconstructions is when the norm of the residual
operator reaches the noise norm (known in the regularization
literature as the discrepancy principle
\cite{vogel2002computational}).
\begin{figure}[t]
\hspace{-1.13cm}
\begin{tabular}{cc}
\epsfig{file=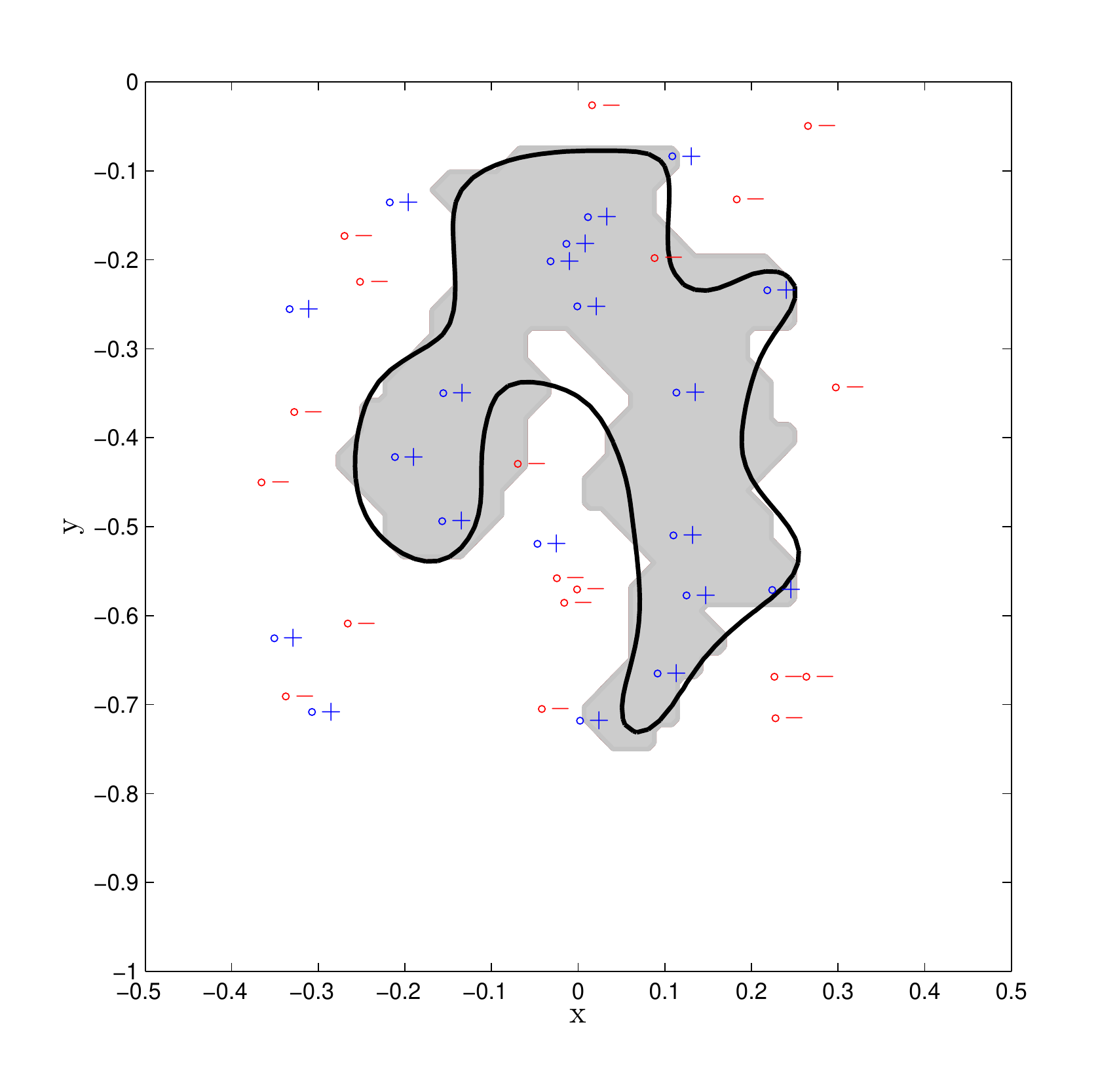,width=0.4\linewidth,clip=a}\hspace{-0.6cm}
\epsfig{file=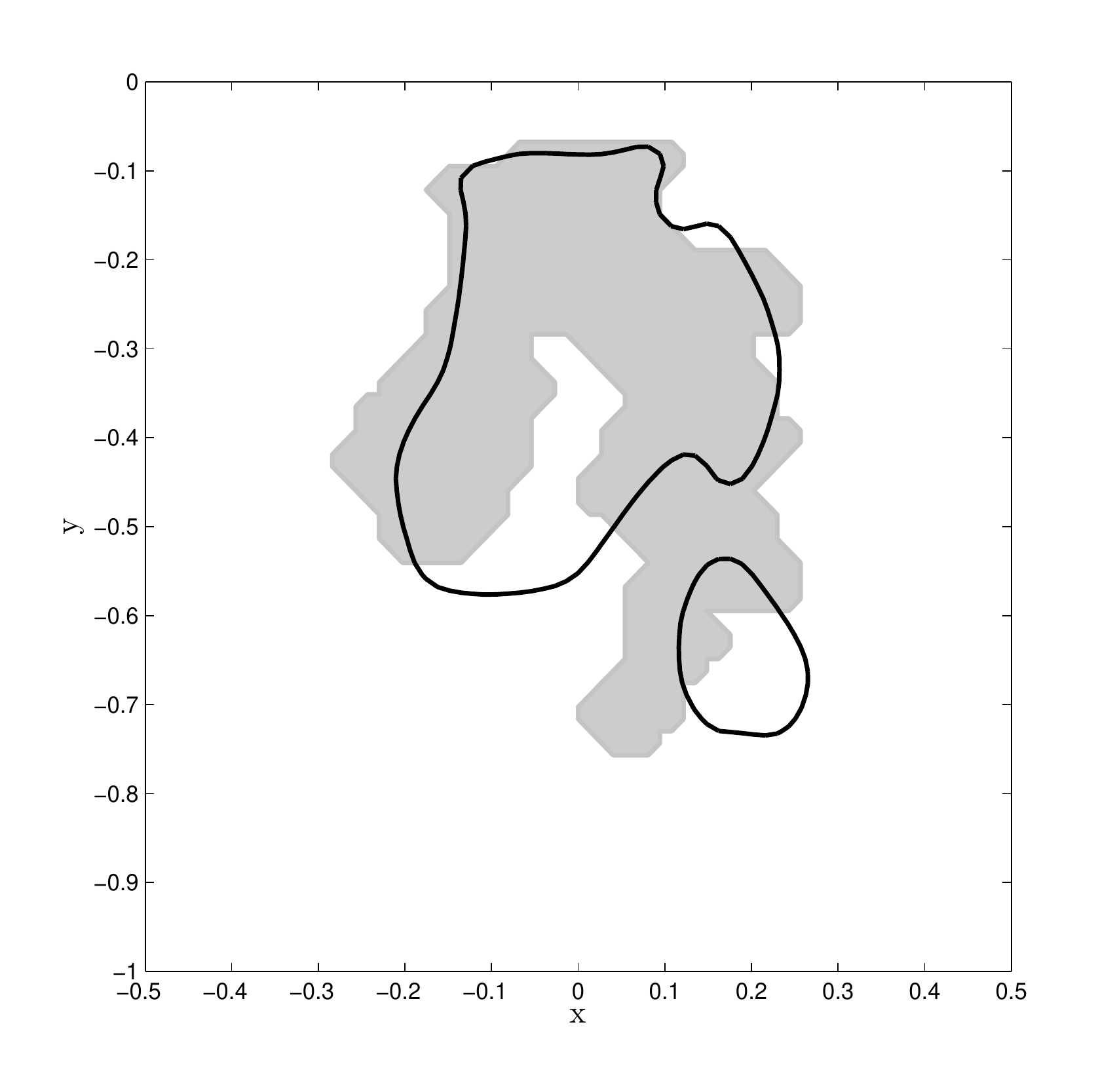,width=0.4\linewidth,clip=a}\hspace{-0.6cm}
\epsfig{file=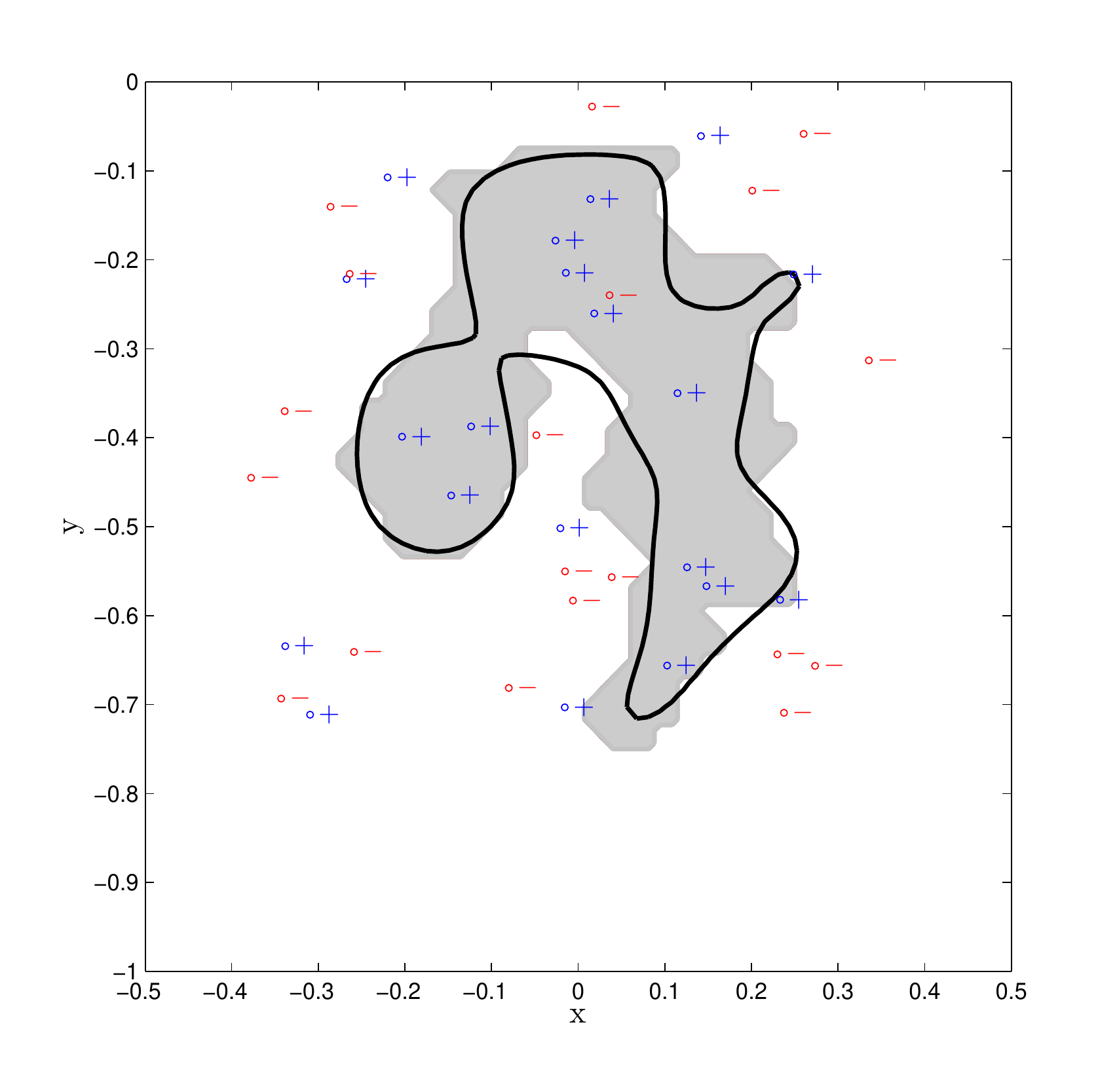,width=0.4\linewidth,clip=a}
\end{tabular}
\caption{(a) Left: The result of reconstructing the shape using
PaLS approach after 26 iterations, the centers
$\boldsymbol{\chi}_j$ and their corresponding weight signs are
shown. (b) Center: Result of using the traditional level set
method after 39 iterations (c) Right: The result of reconstructing
both the shape and the binary conductivity values after 32
iterations}\label{figex2}
\end{figure}
\begin{figure}[h]
\begin{tabular}{cc}
\epsfig{file=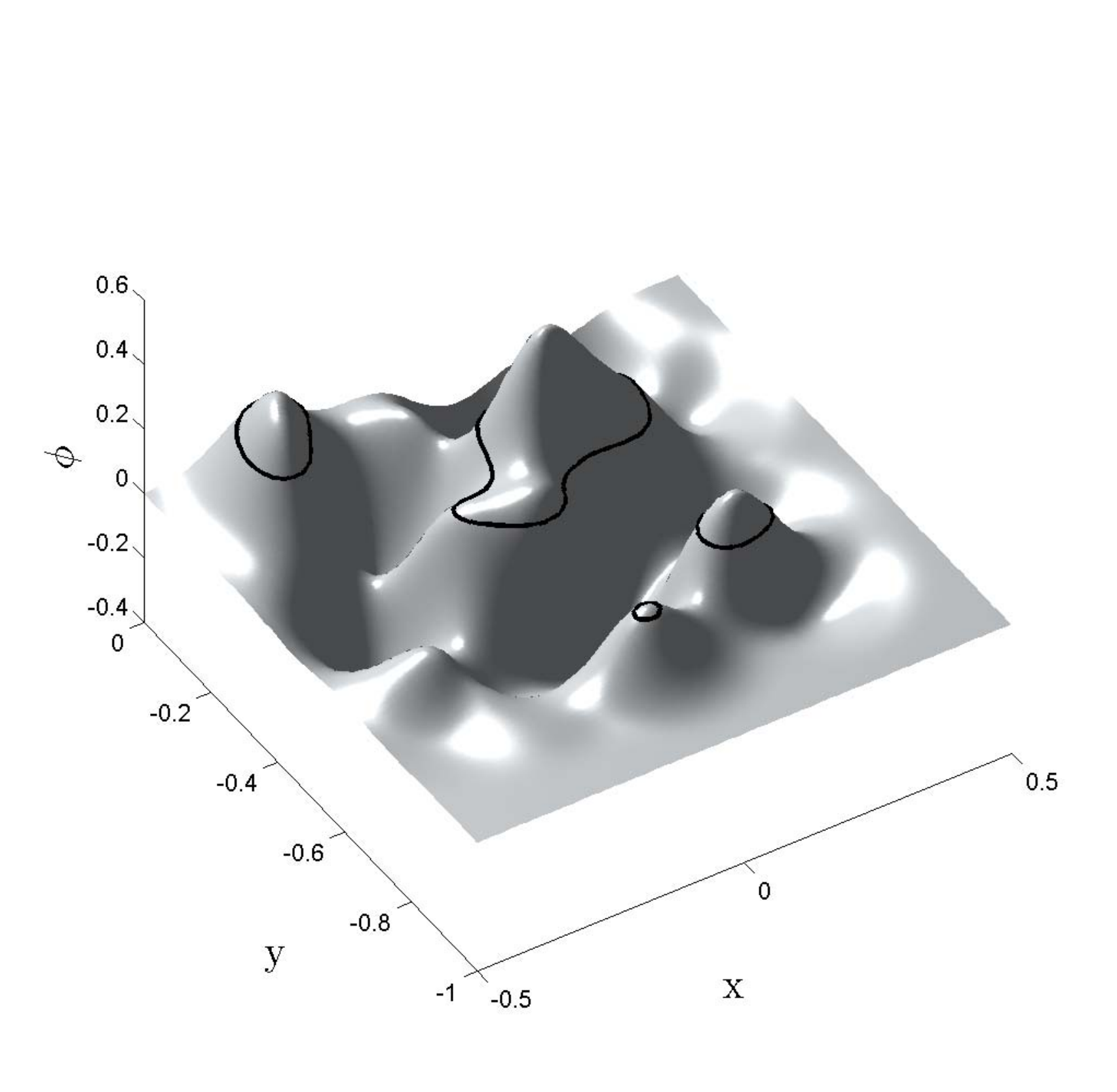,width=0.5\linewidth,clip=a}
\epsfig{file=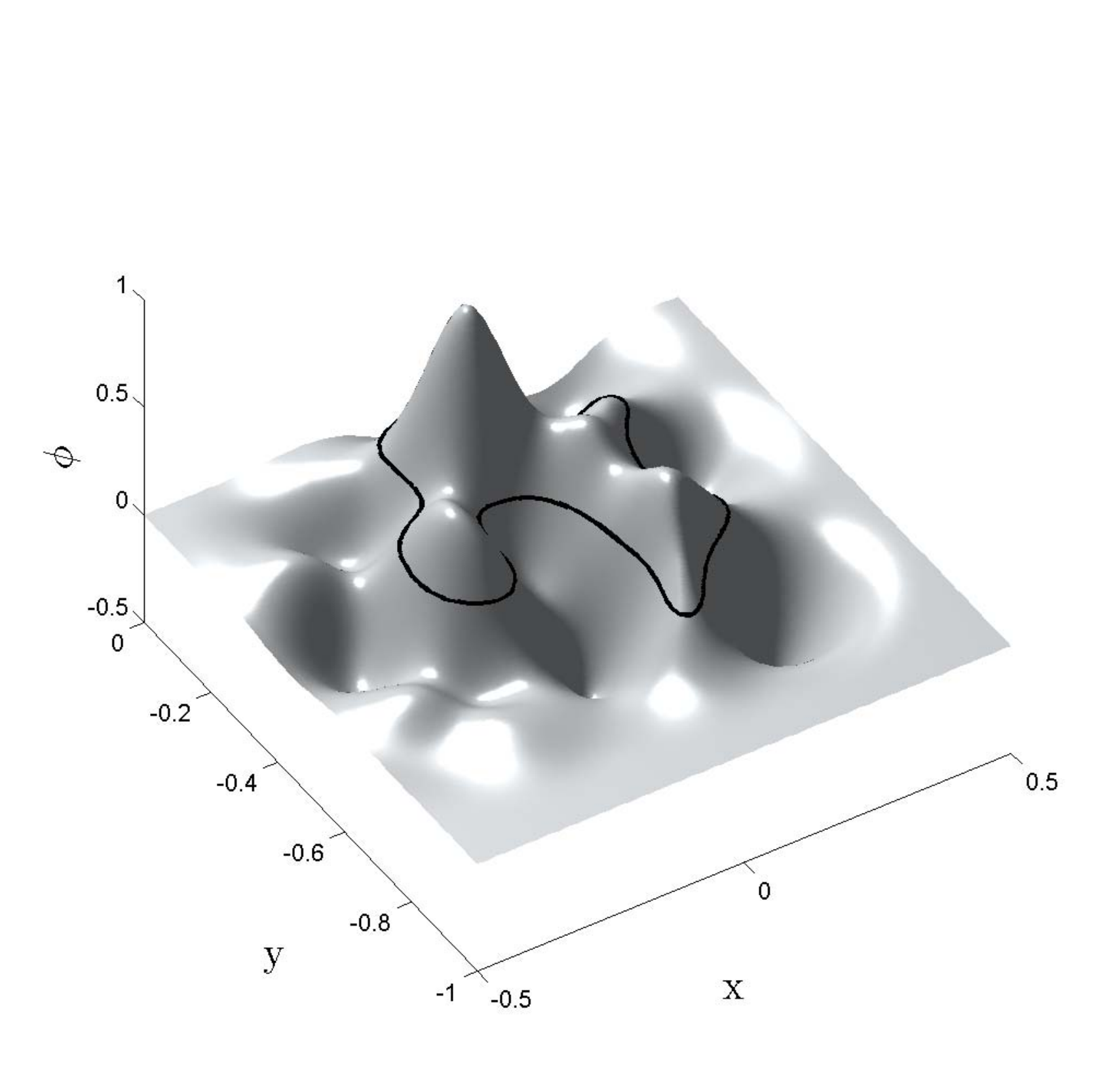,width=0.5\linewidth,clip=a}
\end{tabular}
\caption{(a) Left: The initial PaLS function (b) Right: The final
PaLS function for the shape-only reconstruction of Figure
\ref{figex2}.a}\label{figex3}
\end{figure}
\begin{figure}[h]
\centering \epsfig{file=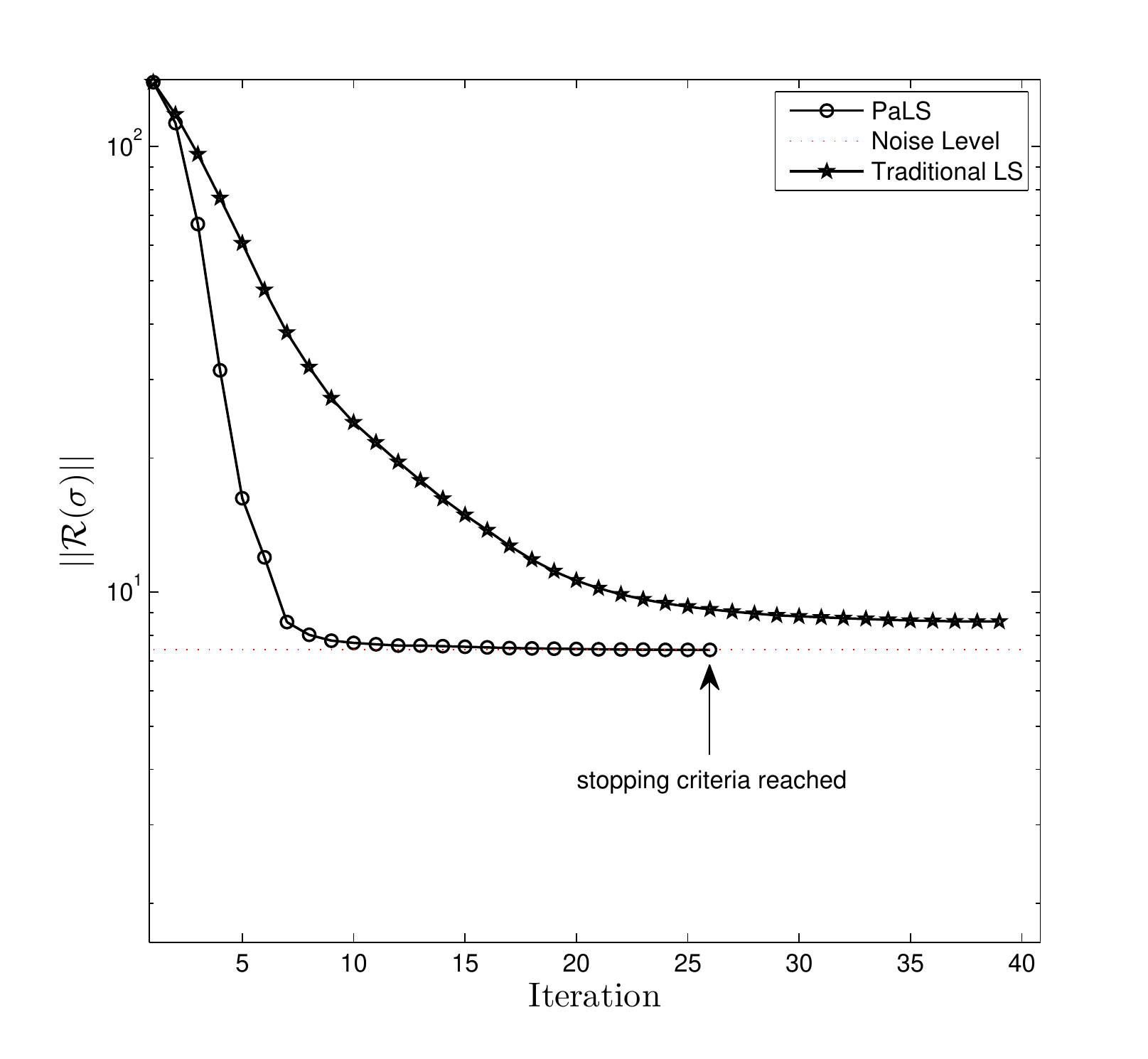,width=0.6\linewidth,clip=a}
\caption{Residual error reduction through the iterative process
for the shape only reconstruction, using PaLS approach and the
traditional level set method }\label{figex4}
\end{figure}

Figure \ref{figex2}.a shows the shape reconstruction using the
PaLS approach, and assuming the values $\sigma_i$ and $\sigma_o$
are a priori known. Figure \ref{figex2}.b shows the result of the
same problem using a typical Gauss-Newton approach with the level
set function defined as a signed distance function over the pixels
and a smoothing regularization added as explained in
\cite{soleimani2006level}. This algorithm only considers shape
reconstruction (i.e., $\sigma_i$ and $\sigma_o$ are considered
known) and it is initialized with the same contour as the initial
$c$-level set contour of the PaLS approach. As the results in
Figure \ref{figex2}.a and \ref{figex2}.b show, the PaLS approach
performs well in reconstructing major shape characteristics, while
the traditional level set approach fails to provide a good
reconstruction in low sensitivity regions close to the bottom and
does not capture the concavity. Figure \ref{figex2}.c shows the
result of reconstructing both the shape and the anomaly values
using the PaLS approach, where this time the PaLS evolution takes
slightly more iterations (32 iterations verses 26), but the
resulting reconstruction still well represents the shape. The
initial values used for the conductivity values are
$\sigma_i^{(0)}=0.01\;\mbox{Sm}^{-1}$ and
$\sigma_o^{(0)}=0.005\;\mbox{Sm}^{-1}$ and the final resulting
values are $\sigma_i^{(32)}=0.056\;\mbox{Sm}^{-1}$ and
$\sigma_o^{(32)}=0.010\;\mbox{Sm}^{-1}$, which show a good match
with the real values. To illustrate the behavior of the PaLS
function, in Figure \ref{figex3} we have shown the initial and
final PaLS functions for the shape only reconstruction. Also to
compare the convergence behaviors of the PaLS approach and the
traditional level set approach, in Figure \ref{figex4} we show the
residual error through the evolution steps for both methods. Using
the PaLS approach the stopping criteria is met after 26 iterations
while traditional level set method reaches a local minima after 39
iterations (the updates after 39 iterations become so small that
it stops evolving further).

\subsection{X-ray Computed Tomography}
As the second example and a mildly ill-posed problem, we consider
X-ray Computed Tomography (CT) \cite{davison1983ill}. In this
imaging technique, X-ray photons are transmitted through the
object of interest and the intensity of transmitted ray is
measured at the boundaries to reconstruct the object's attenuation
coefficient. The contrast between the attenuation characteristics
of different materials can provide structural information about
the object being imaged. X-ray CT is among the most well known
methods for imaging the body tissue in medical applications
\cite{kalender2006x}.

For an X-ray beam transmitting along a line $\mathcal{L}_k$ in the
tissue, the photon intensity $\mathcal{X}_k$ measured at the
detector side of the line can be written as
\begin{equation}
\label{eq2-1} \mathcal{X}_k=\int I_k(\mathcal{E})\exp
\big(-\int_{\mathcal{L}_k}\alpha(\mathbf{x},\mathcal{E})\;\mbox{d}
\mathbf{x}\big)\;\;\mbox{d} \mathcal{E}
\end{equation}
where $\alpha(\mathbf{x},\mathcal{E})$ denotes the attenuation
coefficient, in general as a function of the position $\mathbf{x}$
and the energy of the incident ray $\mathcal{E}$, and
$I_k(\mathcal{E})$ denotes the incident ray energy spectrum. In
case of a monoenergetic beam as
$I_k(\mathcal{E})=I_{0,k}\delta(\mathcal{E}-\mathcal{E}_0)$, a
measured quantity related to the photon intensity may be defined
as
\begin{equation}
\label{eq2-2}
u_k:=-\log(\frac{\mathcal{X}_k}{I_{0,k}})=\int_{\mathcal{L}_k}\alpha(\mathbf{x})
\mbox{d} \mathbf{x}.
\end{equation}
Equation (\ref{eq2-2}) simply relates the measurements to the
Radon transform of the attenuation coefficient
$\alpha(\mathbf{x})$ in a monoenergetic scenario. The quantities
$u_k$ are actually what is considered as the data in CT imaging.

The Fr\'{e}chet derivative of the CT measurements with respect to
the attenuation coefficient is expressed as
\begin{equation}
\label{eq2-3} \frac{\mbox{d}u_{k}}{\mbox{d}\alpha}\;
[\delta\alpha]=\int_{\mathcal{L}_k}\delta
\alpha\;\mbox{d}\mathbf{x}.
\end{equation}
Considering $\mathbf{u}$ to be the set of CT data collected along
different paths $\mathcal{L}_k$, for $k=1,2,\cdots N$, and
$\mathcal{M}(\alpha)$ as the forward model mapping the attenuation
to the CT data set, based on (\ref{eq2-3}) the sensitivity of the
residual operator
$\mathcal{R}(\alpha)=\mathcal{M}(\alpha)-\mathbf{u}$ with respect
to a perturbation in $\alpha$ can be written as
\begin{equation}
\label{eq2-4} \mathcal{R}'(\alpha)[\delta\alpha] = \left(
\begin{array}{c}
\int_{\mathcal{L}_1} \delta \alpha \;\mbox{d}\mathbf{x}\\ \vdots\\
\int_{\mathcal{L}_N} \delta \alpha \;\mbox{d}\mathbf{x}\\
\end{array} \right),
\end{equation}
which we need for the PALS evolution process.
\begin{figure}[t]
\hspace{-1.3cm}
\begin{tabular}{cc}
\epsfig{file=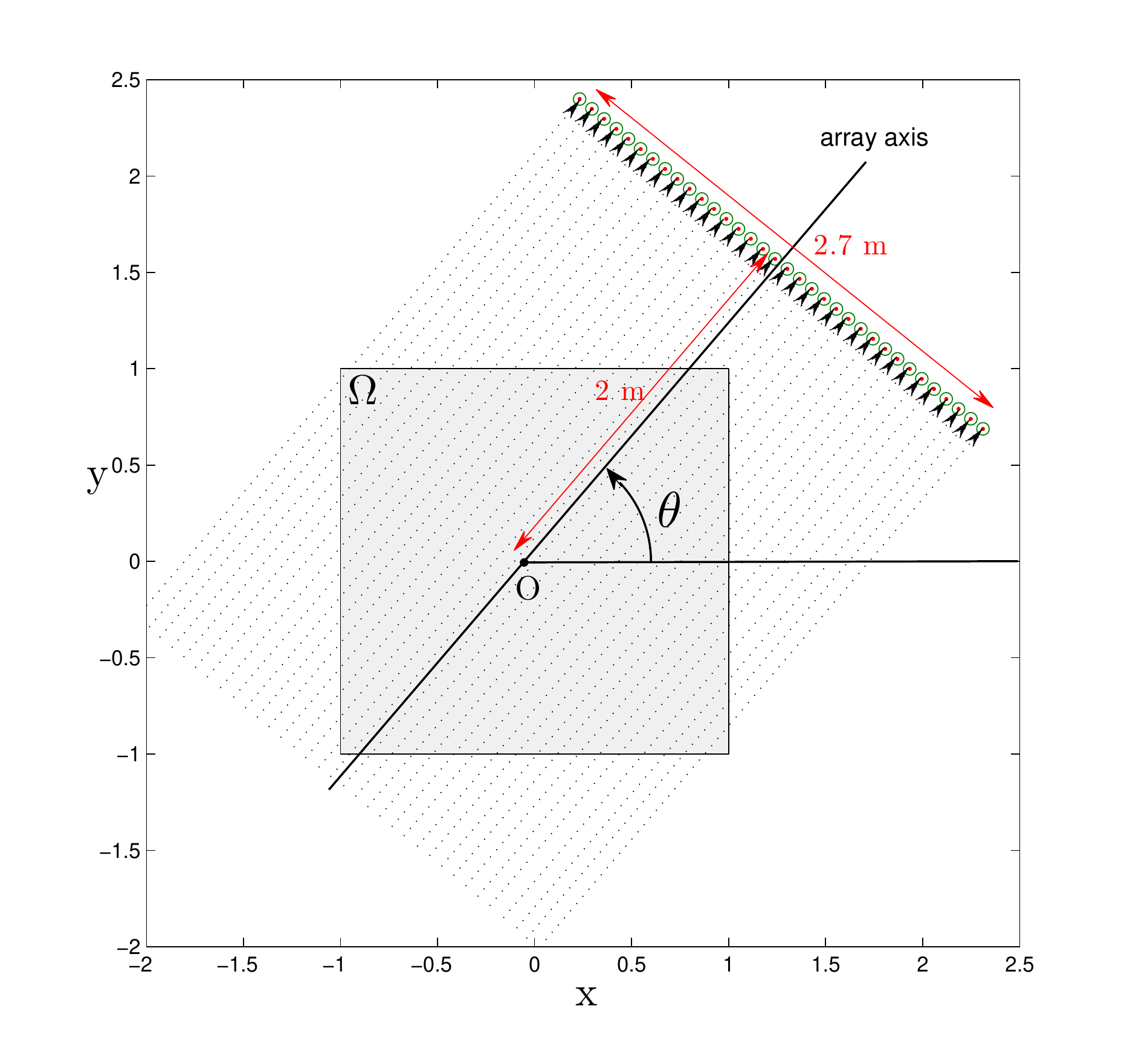,width=0.58\linewidth,clip=a}
\epsfig{file=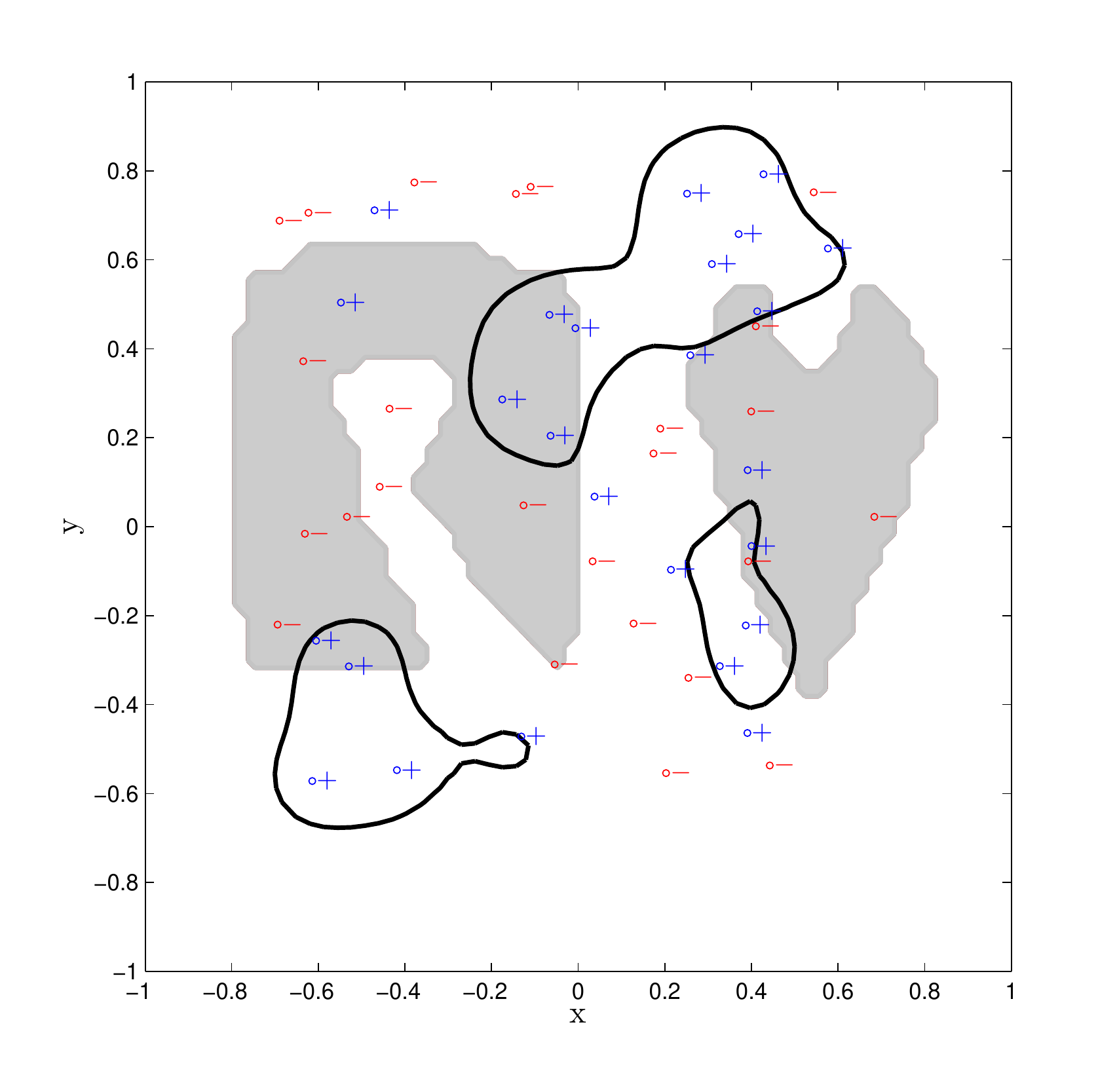,width=0.56\linewidth,clip=a}
\end{tabular}
\caption{(a) Left: The 2D X-ray CT imaging setup (b) Right: The
gray region shows the attenuation shape to be reconstructed. The
dots with ``+" and ``-" signs correspond to the centers of
positive and negative weighted bumps in the initial state of the
problem. The black contour is the resulting $c$-level set of the
initial PaLS function }\label{figex5}
\end{figure}

We consider 2D imaging over a square of $2 \; \mbox{m}\times 2\;
\mbox{m}$, i.e., $\Omega=[-1,1]\times[-1,1]$ in the $x-y$ plane,
as shown in Figure \ref{figex5}.a. The region outside $\Omega$ is
assumed to have zero attenuation. The X-ray beams are emitted as
parallel beams and the measurements are performed through an
equally spaced linear array of 34 sensors, the vertical axis of
which makes an angle $\theta$ with the $x$ axis. For a full view
X-ray CT we have $\theta_0\leq \theta \leq \theta_0+\pi$, while in
a limited view scenario, and hence a more ill-posed problem,
$\theta$ varies in a smaller angular domain. The shape to be
reconstructed is shown in Figure \ref{figex5}.b, which is brought
here from \cite{tai2004survey} as a rather complex shape to
examine the flexibility of the PaLS approach. The region $\Omega$
is discretized into $64\times 64$ uniform grid. The binary
attenuation values are $\alpha_i=2.5 \;\mbox{cm}^{-1}$ and
$\alpha_o=1 \;\mbox{cm}^{-1}$. For the purpose of imaging, the
forward model measurements are performed at every $1$ degree
angle. For the PaLS representation, we use $m_0=50$ bumps with the
centers $\boldsymbol{\chi}_j$, distributed randomly as shown in
Figure \ref{figex5}.b. For this example we use slightly more bumps
due to the better posed nature of the problem (at least in the
full view case) and the more complex shape. Again to roughly
carpet the domain $\Omega$ with the bumps we uniformly take the
dilation factors to be $\beta_j=2.5$. The other PaLS settings used
are similar to those of the ERT example in the previous section.
\begin{figure}[t]
\hspace{-.5cm}
\begin{tabular}{cc}
\epsfig{file=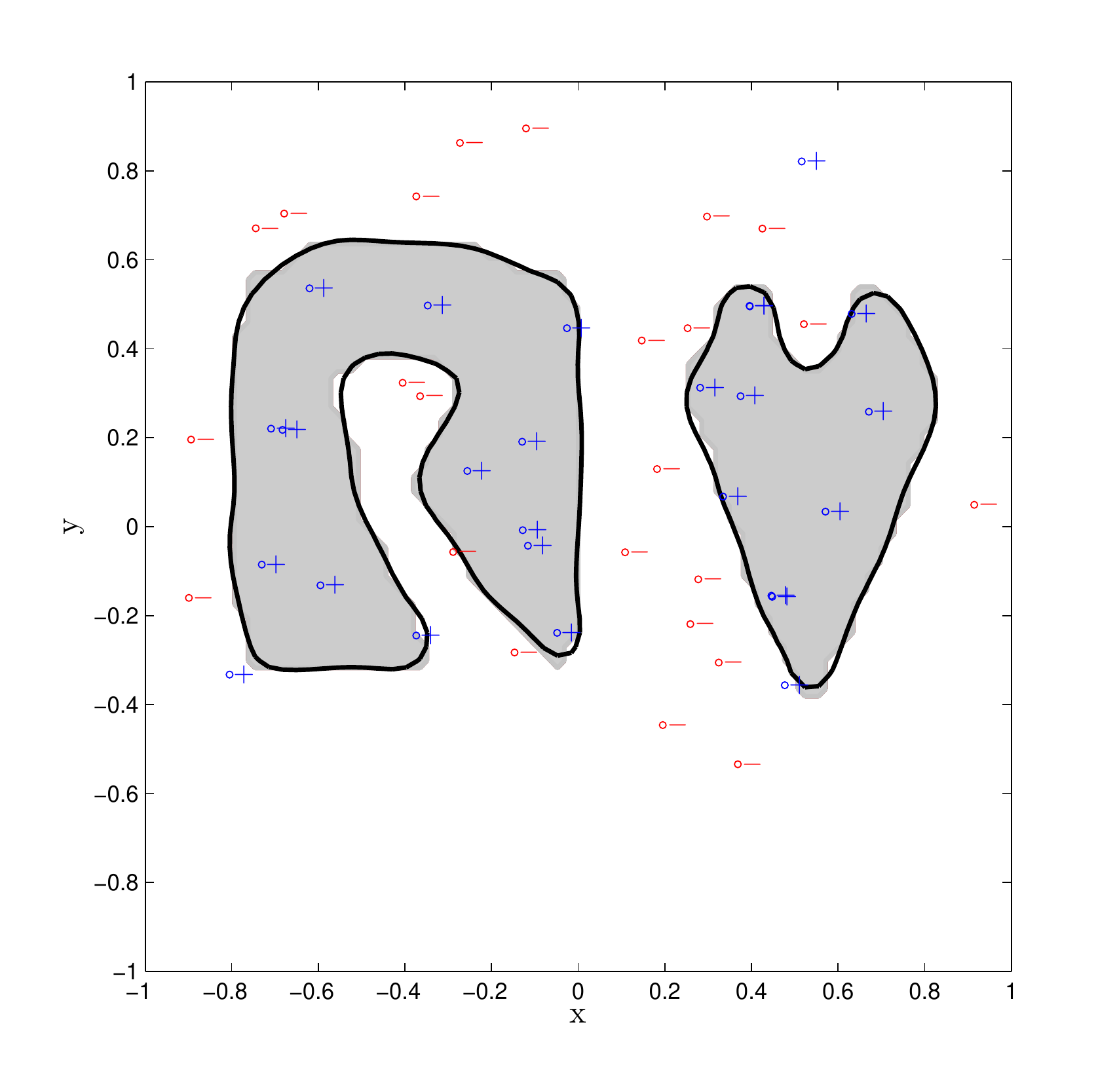,width=0.5\linewidth,clip=a}
\epsfig{file=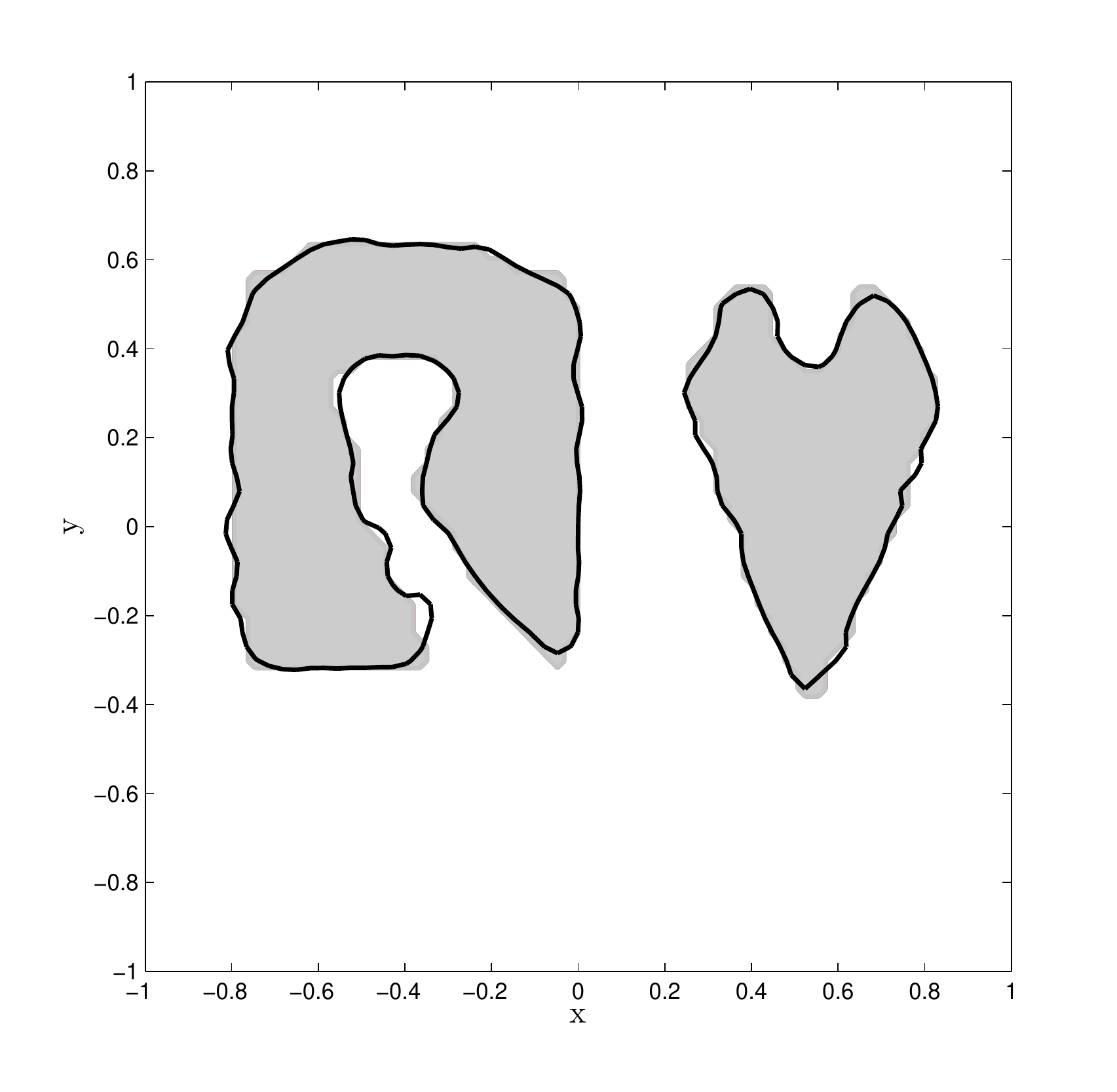,width=0.5\linewidth,clip=a}
\end{tabular}
\caption{(a) Left: The result of reconstructing both the shape and
the anomaly values in a $5\%$ noise case using the PaLS approach,
convergence achieved after 20 iterations (b) Right: Result of
using the same data set with traditional level set method to
recover shape only, assuming attenuation coefficients are known,
convergence achieved after 14 iterations }\label{figex6}
\end{figure}

\begin{figure}[th]
\hspace{-.5cm}
\begin{tabular}{cc}
\epsfig{file=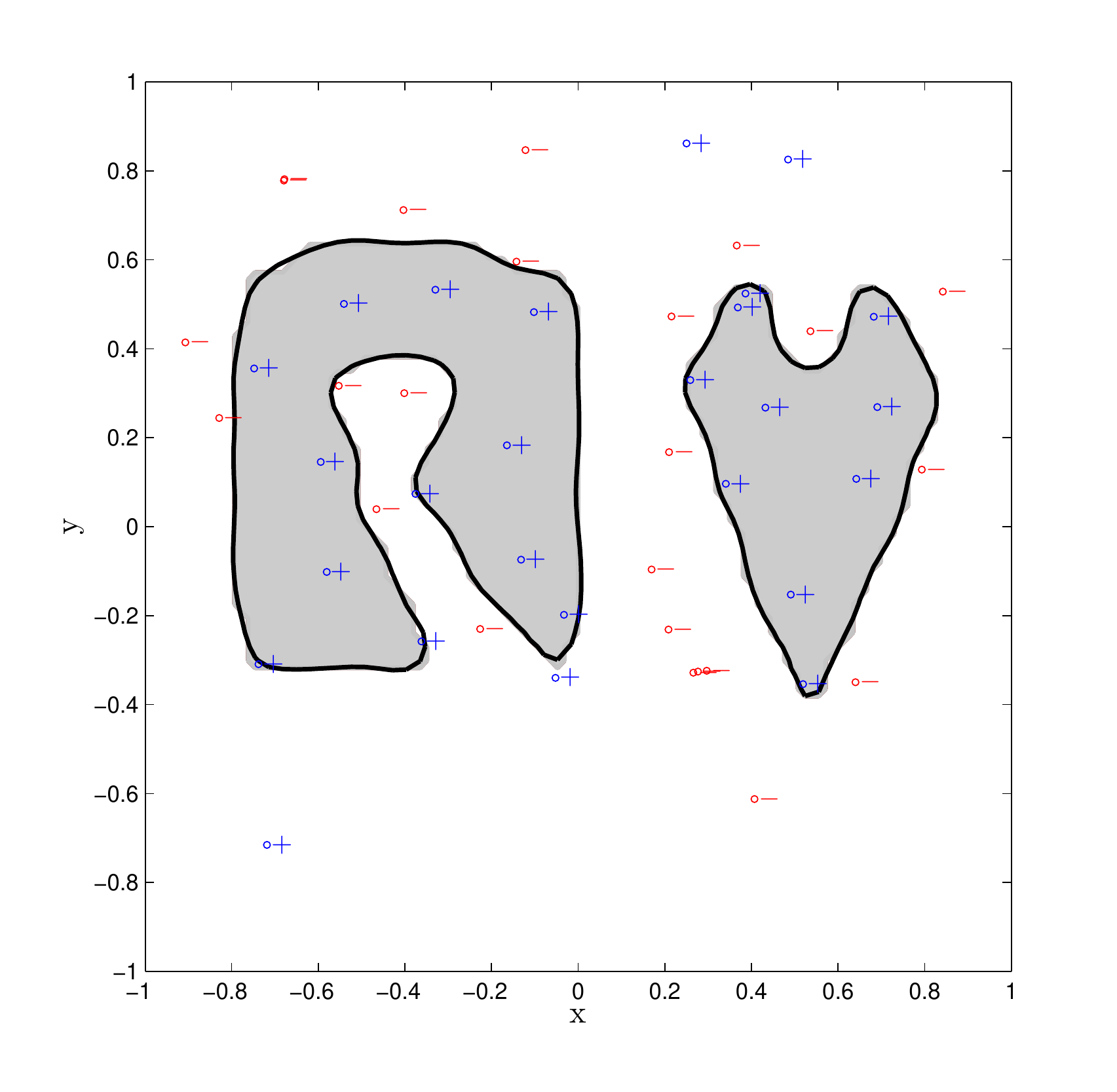,width=0.5\linewidth,clip=a}
\epsfig{file=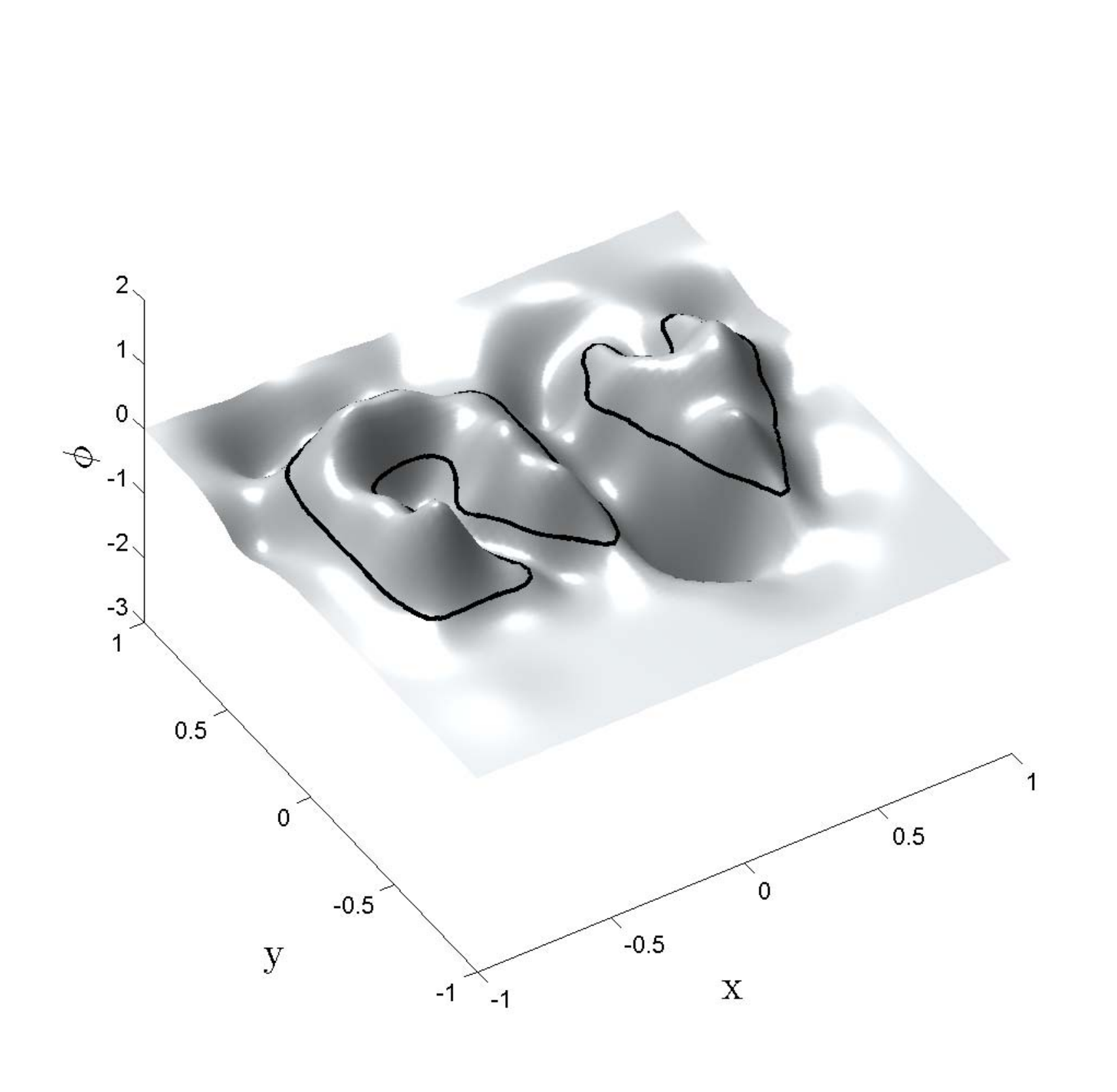,width=0.5\linewidth,clip=a}
\end{tabular}
\caption{(a) Left: The result of PaLS shape reconstruction in a
$1\%$ noise case to show the flexibility of following shape
details due to the pseudo-logical property, convergence achieved
after 42 iterations (b) Right: The corresponding PaLS function
}\label{figex7}
\end{figure}

\begin{figure}[h]
\hspace{-.5cm}
\begin{tabular}{cc}
\epsfig{file=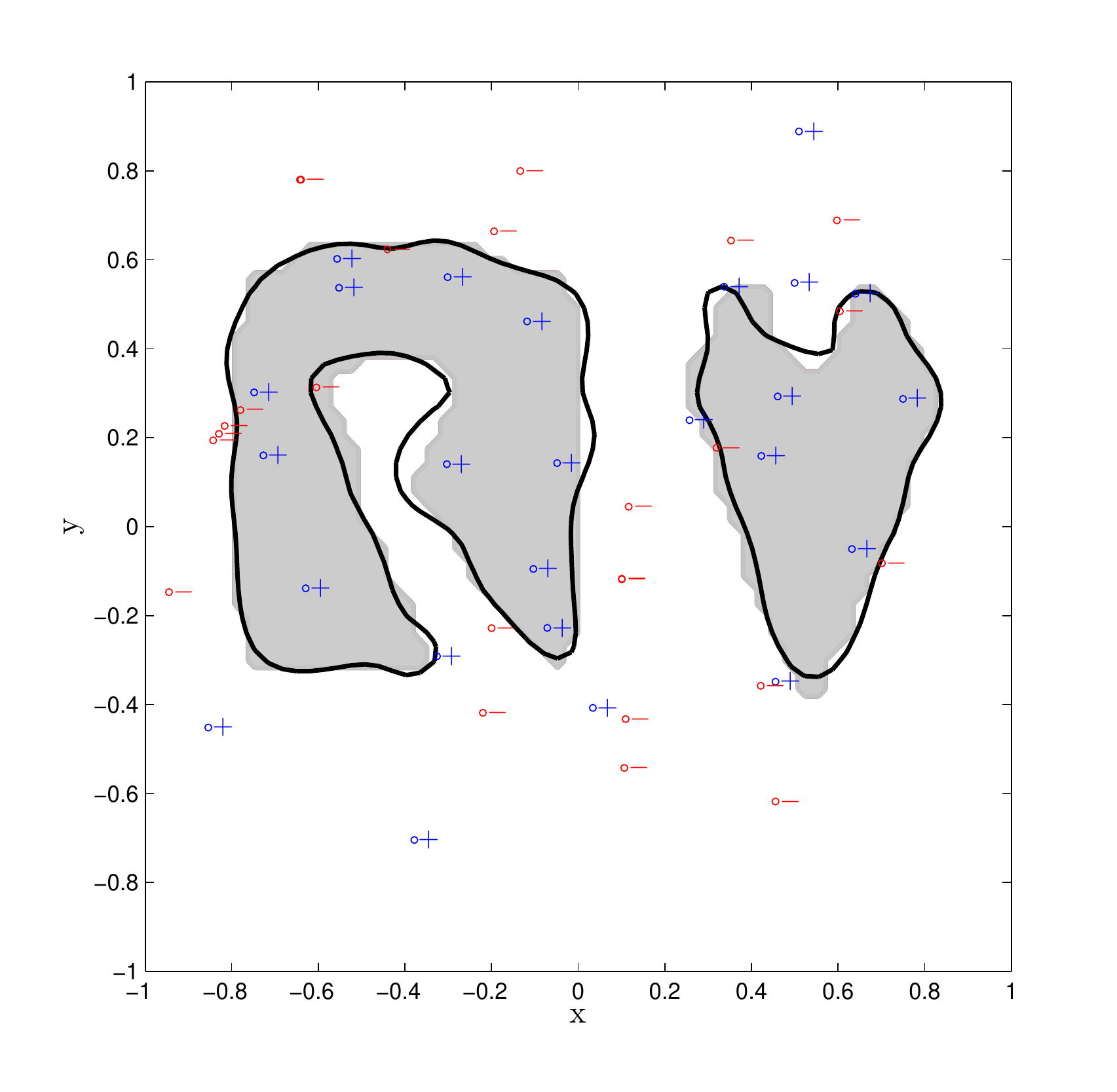,width=0.5\linewidth,clip=a}
\epsfig{file=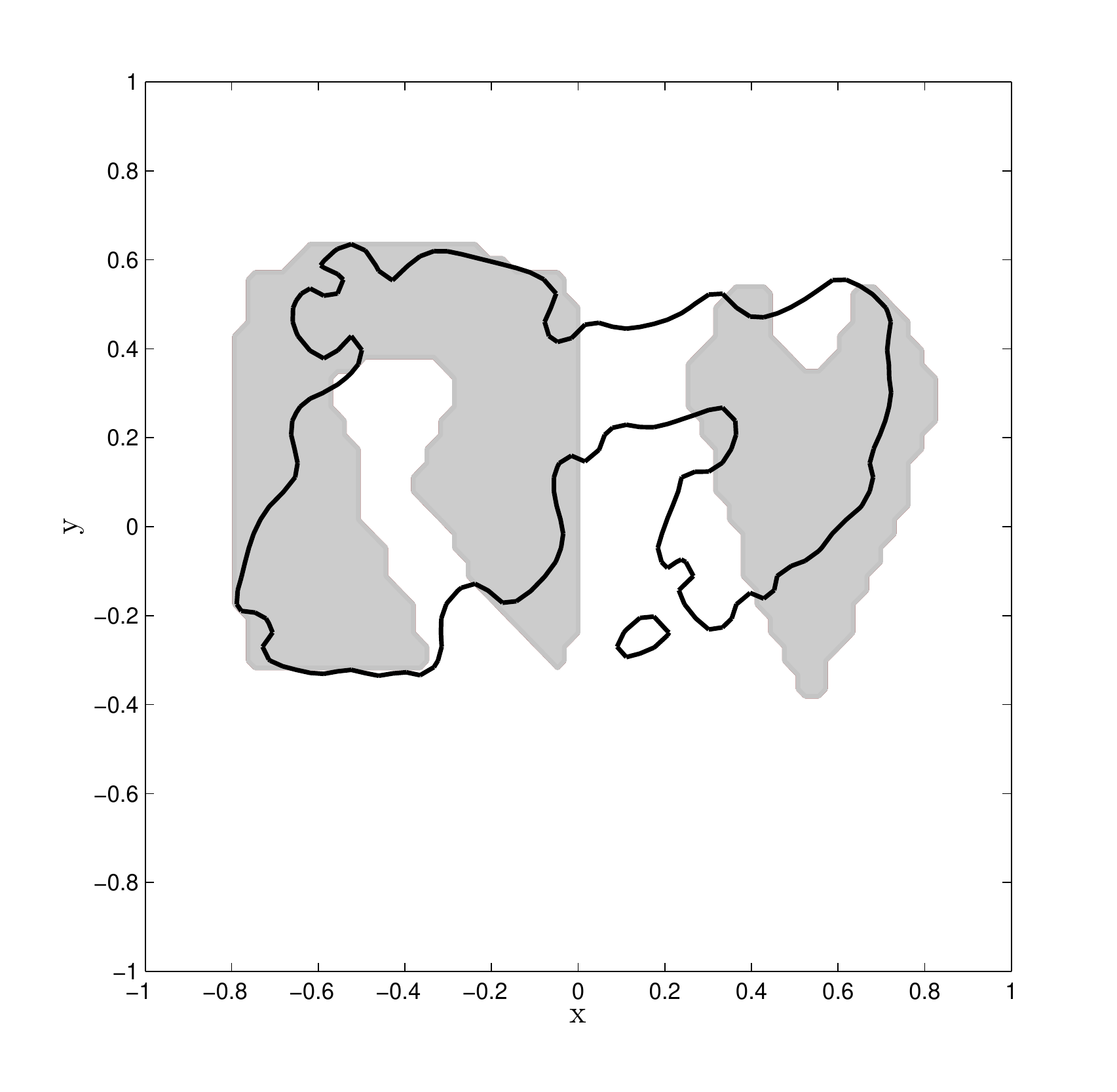,width=0.5\linewidth,clip=a}
\end{tabular}
\caption{(a) Left: The PaLS shape reconstruction in a limited view
X-ray CT, convergence achieved after 49 iterations (b) Right: The
performance of the traditional level set method for the same data
set, reaching a local minima after 32 iterations}\label{figex8}
\end{figure}

Figure \ref{figex6}.a shows the result of reconstructing both the
shape and the attenuation values for a full view experiment where
$0<\theta<\pi$. The synthetic data are obtained from the true
attenuation map and by adding $5\%$ Gaussian noise to the
measurements. The initial values used for the attenuations are
$\alpha_i^{(0)}=1.5\;\mbox{cm}^{-1}$ and
$\alpha_o^{(0)}=0.5\;\mbox{cm}^{-1}$ and the discrepancy principle
stopping criteria is met after 20 iterations resulting the
corresponding shape and the final values of
$\alpha_i^{(20)}=2.494\;\mbox{cm}^{-1}$ and
$\alpha_o^{(20)}=0.997\;\mbox{cm}^{-1}$, which are very well
matched with the real quantities. We should mention here that when
the attenuation values are assumed to be known, the number of
iterations for shape reconstruction is only 12. The same data set
is used to reconstruct the shape only, using the traditional level
set method used in the ERT example and initialized with the same
contour as the $c$-level set contour shown in Figure
\ref{figex5}.b. The result of this reconstruction is shown in
Figure \ref{figex6}.b. Due to the better posed nature of this
problem, the pixel based level set method also provides a
successful reconstruction for the geometry only. Note, however
that the PaLS approach is solving a more challenging problem
(i.e., reconstructing the attenuation values as well as the
shape), thanks to the pseudo-logical behavior of the bumps, the
resulting shape follows the same detail levels as a pixel based
level set function. This shaping capability of the PaLS is more
highlighted when we use a less noisy data ($1\%$ Gaussian noise).
The result is shown in Figure \ref{figex7} where the convergence
is achieved after 42 iterations, and the resulting contours
follows the true shape in high level of details. Finally as a more
challenging and very ill-posed problem, we now consider a limited
view scenario where the angular domain is limited to
$\pi/4<\theta<3\pi/4$, and fewer rays cross the anomalies. Two
percent additive Gaussian noise is added to the synthetic data
obtained from the true attenuation map. With the same problem
setting as the full view case, Figure \ref{figex8}.a shows the
result of reconstructing the shape using the PaLS approach.
However, as shown in Figure \ref{figex8}.b, the traditional level
set method applied to this problem given the attenuation values,
fails to provide a complete reconstruction and stops further shape
enhancement after reaching a local minima.

\subsection{Diffuse Optical Tomography}

As the third application we consider diffuse optical tomography
(DOT). In this imaging method data are obtained by transmitting
near-infrared light into a highly absorbing and scattering medium
and then recording the transmitted light. As the inverse problem
of interest the photon flux is measured on the surface to recover
the optical properties of the medium such as absorption and/or
reduced scattering functions. A well known application of this
method is breast tissue imaging, where the differences in the
absorption and scattering may indicate the presence of a tumor or
other anomaly \cite{boas2001imaging}. For given absorption and
scattering functions a number of mathematical models have been
proposed in the literature to determine the synthetic data (photon
flux) \cite{arridge}. We focus on the frequency-domain diffusion
model in which the data is a non-linear function of the absorption
and scattering functions.

Consider $\Omega$, a square with a limited number of sources on
the top and a limited number of detectors on either the top or
bottom or both. We use the diffusion model in \cite{arridge}
modified for the 2D case where the photon flux
$u(\mathbf{x};\omega)$ is related to input $s(\mathbf{x};\omega)$
through
\begin{equation}
\label{eqe3-1} -\nabla\cdot \beta(\mathbf{x}) \nabla
u(\mathbf{x};\omega) + \alpha(\mathbf{x}) u(\mathbf{x};\omega)
+i\frac{\omega}{v}u(\mathbf{x};\omega) = s(\mathbf{x};\omega) ,
\end{equation}
with the Robin boundary conditions
\begin{equation}
\label{eqe3-2}u(\mathbf{x};\omega) + 2\beta(\mathbf{x})
\frac{\partial
    u(\mathbf{x};\omega)}{\partial \nu} = 0, \qquad \mathbf{x} \in
    \partial \Omega.
\end{equation}
Here, $\beta(\mathbf{x})$ denotes the diffusion, which is related
to the reduced scattering function $\mu_s'(\mathbf{x})$, by
$\beta(\mathbf{x}) = 1/\big(3 \mu_s'(\mathbf{x})\big)$,
$\alpha(\mathbf{x})$ denotes absorption and $\partial/\partial
\nu$ denotes the normal derivative. We also have $i=\sqrt{-1}$,
$\omega$ as the frequency modulation of light and $v$ being the
speed of light in the medium. Knowing the source and the functions
$\alpha(\mathbf{x})$ and $\beta(\mathbf{x})$, we can compute the
corresponding $u(\mathbf{x};\omega)$ everywhere, in particular, at
the detectors.

As the inverse problem, we consider a case that the reduced
scattering function is known and we want to reconstruct the
absorption $\alpha(\mathbf{x})$ from the data. Again for this
problem we consider a point source
$s(\mathbf{x};\omega)=\delta(\mathbf{x}-\mathbf{x}_s)$ for
$\mathbf{x}_s\in\partial\Omega$, which results in a photo flux
$u_s(\mathbf{x};\omega)$ in $\Omega$. For a measurement at
$\mathbf{x}_d\in \Omega$ as
\begin{equation}
\label{eqe3-3}u_{ds}=\int_\Omega
u_s(\mathbf{x};\omega)\delta(\mathbf{x}-\mathbf{x}_d)\mbox{d}\mathbf{x},
\end{equation}
the variations with respect to the variations in the absorption
can be written as \cite{arridge}
\begin{equation}
\label{eqe3-4}\frac{\mbox{d}u_{ds}}{\mbox{d}\alpha}[\delta
\alpha]=\int_\Omega \delta\alpha
\;u_s(\mathbf{x};\omega)\;u_d(\mathbf{x};\omega)\mbox{d}\mathbf{x},
\end{equation}
with $u_d$ being the adjoint field caused by having the point
source at $\mathbf{x}_d$. Identical to the notation used for the
ERT problem, consider
$\mathcal{R}(\alpha)=\mathcal{M}(\alpha)-\mathbf{u}$ as the
residual operator where the data vector $\mathbf{u}$ contains the
measurements of $N_\ell$ detectors from $M$ experiments. Based on
(\ref{eqe3-4}) the $\ell^{\mbox{th}}$ block of
$\mathcal{R}'(\alpha)[.]$ corresponding to the $\ell^{\mbox{th}}$
experiment is
\begin{equation}
\label{eq3-5} \mathcal{R}'_\ell(\alpha)[\delta\alpha] = \left(
\begin{array}{c}
\int_\Omega \delta \alpha \; u_\ell \;
 u_{1}^{\ell} \; \mbox{d}\mathbf{x}\\
\vdots\\ \int_\Omega \delta \alpha \; u_\ell \;  u_{N_\ell}^{\ell}
\; \mbox{d}\mathbf{x}
\end{array} \right),
\end{equation}
with $u_\ell$ denoting the photon flux in the $\ell^{\mbox{th}}$
experiment and $u_i^\ell$ being the adjoint field.
\begin{figure}[th]
\centering \epsfig{file=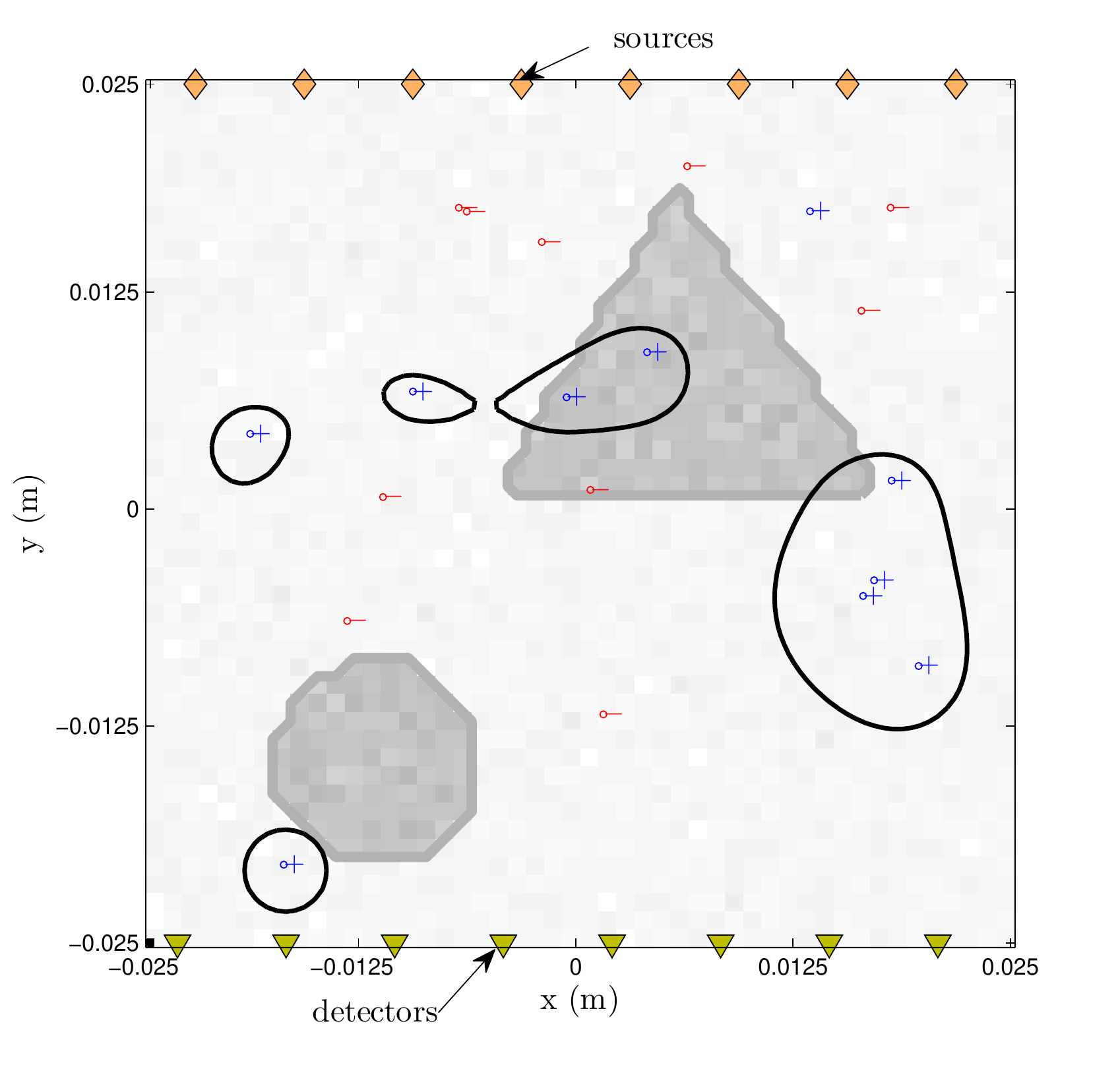,width=0.5\linewidth,clip=a} \caption{The
2D DOT imaging setup with the sources and detectors setting. The
gray region shows the absorption shape, while $2\%$ Gaussian noise
is added to the absorption as a heterogeneity. The dots with ``+"
and ``-" signs correspond to the centers of positive and negative
weighted bumps in the initial state of the problem. The black
contour is the resulting $c$-level set of the initial PaLS
function }\label{figex9}
\end{figure}
\begin{figure}[h]
\hspace{-.5cm}
\begin{tabular}{cc}
\epsfig{file=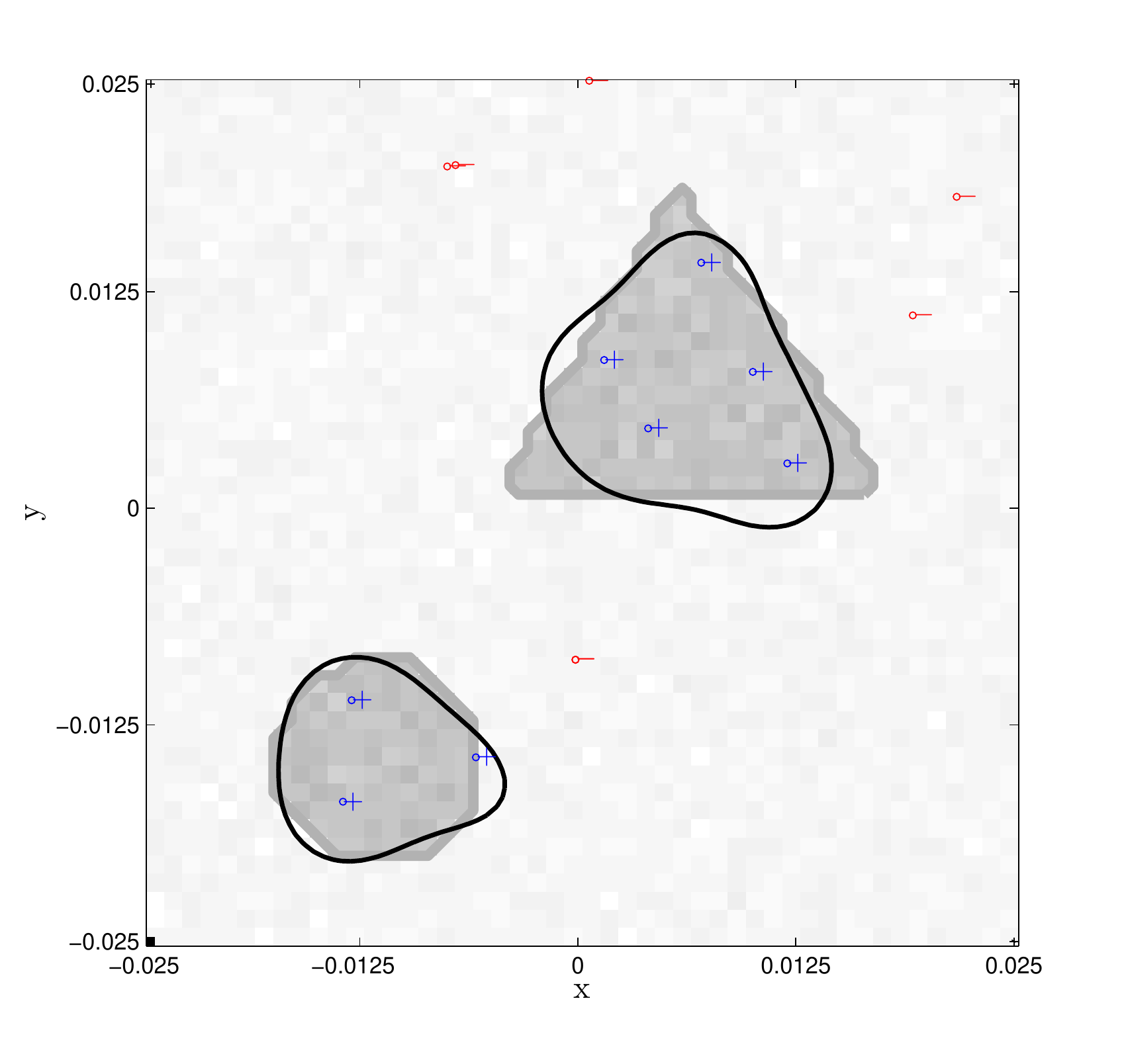,width=0.5\linewidth,clip=a}
\epsfig{file=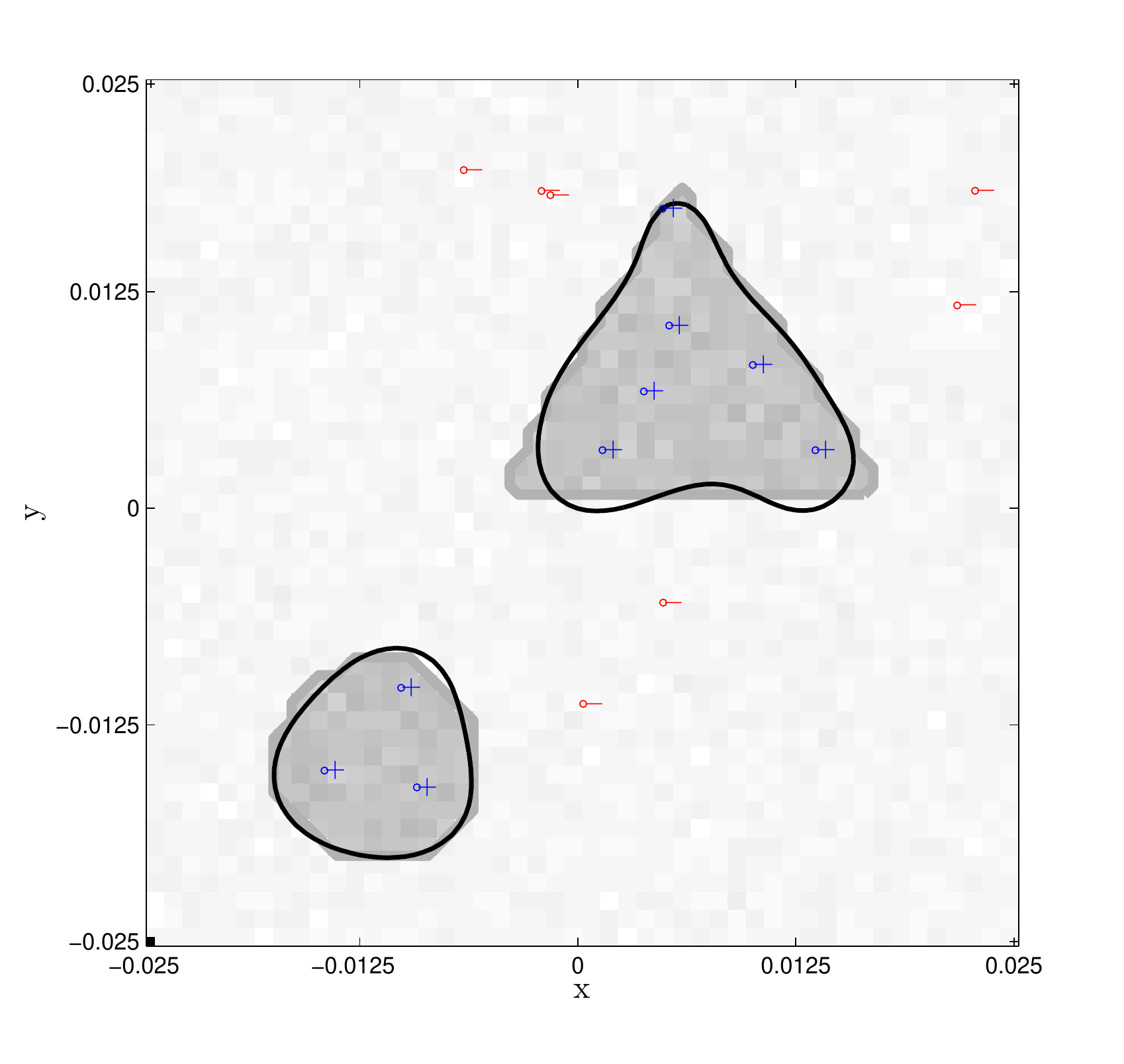,width=0.5\linewidth,clip=a}
\end{tabular}
\caption{(a) Left: The PaLS shape-only reconstruction, stopping
criteria achieved after 152 iterations in case $0.1\%$ data noise
and $2\%$ heterogeneity (b) Right: The results of shape
reconstruction in case of only heterogeneity in the absorption.
The evolution is manually stopped after 200 iterations, showing
the results}\label{figex10}
\end{figure}
For the purpose of this example we consider a very ill-posed
problem where $\Omega$ corresponds to a $5 \mbox{cm}\times 5
\mbox{cm}$ imaging region, and there are only 8 point sources on
the top and 8 detectors at the bottom. In every experiment only
one of the sources is on and the measurements are performed at the
bottom detectors. We have $\mu_s'(\mathbf{x})=6 \;\mbox{cm}^{-1}$
throughout $\Omega$. The measurements are performed at DC mode,
$25$ MHz and $50$ MHz. The absorption binary distribution is shown
in Figure \ref{figex9} with values $\alpha_i=0.015 \;
\mbox{cm}^{-1}$ and $\alpha_o=0.005 \; \mbox{cm}^{-1}$. To make
the problem more challenging we made the absorption distribution
slightly heterogeneous by adding $2\%$ white Gaussian noise to it.
For the PaLS setting, due to the very ill-posed nature of the
problem we use relatively fewer bumps, i.e., $m_0=20$. Similar to
the previous two examples the centers $\boldsymbol{\chi}_j$
corresponding to positive and negative weighted bumps are taken
randomly inside $\Omega$ as shown in Figure \ref{figex9}. Again to
have reasonable initial support radius for the bumps the dilation
factors are set to be $\beta_j=80$. Other PaLS settings are
similar to the previous examples. The forward model is solved
using finite difference method by discretizing $\Omega$ to
$50\times 50$ grid points. Figure \ref{figex10}.a shows the result
of our shape-only reconstruction after 152 iterations, when the
true absorption map is used to generate the data and $0.1\%$
Gaussian noise is added to it. In case of not having any noise in
the data (but still considering the heterogeneity in the
absorption), the results after 200 iterations are shown in figure
\ref{figex10}.b. Although the current problem settings make it
very ill-posed and extremely hard for shape reconstruction, the
PaLS approach shows reasonable reconstructions, providing
important details. For this very ill-posed problem, using the
traditional level set method we failed to generate useful
reconstructions.

\section{Concluding Remarks}\label{sec7}

In this paper we proposed a parametric level set approach to shape
based inverse problems. The basic formulation of the problem is
kept general, emphasizing the fact that numerically, an
appropriate parametrization of the level set function is capable
of reducing the problem dimensionality and intrinsically
regularizes the problem. Through such modelling the level set
function is more likely to remain well behaved, and therefore
there is no need to re-initialize as in traditional level set
methods. Moreover, based on the fact that the number of underlying
parameters in a parametric approach is usually much less than the
number of pixels (voxels) resulting from the discretization of the
level set function, we make a use of a Newton type method to solve
the underlying optimization problem. We specifically considered
using compactly supported radial basis functions, which used as
proposed provides two main advantages besides being parametric;
first, the pseudo-logical property allows for recovery of wide
range of shapes; and second, implicit narrow-banding for the
evolution.

Although in this paper, we only considered compactly supported
radial basis functions as the bump functions with circular
support, a combination of sufficiently smooth bumps with various
types of supports may be considered in applications where there is
more prior information about the shape and its general geometric
structure. An efficient classification of such functions for
various types of shapes and problems (i.e., forming a basis
dictionary) would be a future direction of research. Although in
the examples presented in this work the number of RBFs where
coarsely chosen, we are currently developing some ideas to
adoptively determine the number of basis elements using a sparse
signal processing approach. Clearly this matter is beyond the
scope of this paper and is planned to be presented as a future
work. Moreover, for sake of simplicity in this paper we only
considered 2D problems, whereas the efficiency of a parametric
approach is more pronounced for 3D shape reconstructions where the
contrast between the number of voxels and the number of PaLS
parameters in a parametric approach is more significant.
Efficiently modelling 3D scenarios via both parametric shape
representation and appropriate texture and heterogeneity models is
an important future direction and a continuation of the current
work.

\section*{Acknowledgments}
This work is funded and supported by NSF, the National Science
Foundation, under grant EAR 0838313.

%
%
\end{document}